\documentclass[12pt]{article}
\usepackage{latexsym}

\newtheorem{lemma}{Lemma}[section]
\newtheorem{theorem}[lemma]{Theorem}
  
\newtheorem{corollary}[lemma]{Corollary}

\newtheorem{proposition}[lemma]{Proposition}

\newtheorem{conjecture}[lemma]{Conjecture}

\newcommand{\comment}[1]{}

\newcommand{\text}[1]{\quad\mbox{#1}\quad}
\def\beq{\begin{equation}}\def\eeq{\end{equation}}
\def\beqn{\begin{eqnarray}}\def\eeqn{\end{eqnarray}}
\def\pont{\hspace{-6pt}{\bf.\ }}

\def\qed{\ifhmode\unskip\nobreak\fi\quad\ifmmode\Box\else$\Box$\fi}
\newcounter{petit}[section]
\def\Nset{\hbox{I\hskip-0.20em I\hskip-0.35em N}}
\def\Qset{\hbox{\hbox{Q\hskip-0.525em\lower-0.097ex
\hbox{\vrule height1.47ex width 0.07em}}\hskip0.50em}}

\newcommand{\petit}[1]{
  \refstepcounter{petit}
  \label{#1}
  (\arabic{petit})
}

\title{The chromatic gap and its extremes}

\author{Andr\'as Gy\'arf\'as\thanks{Research supported in part by OTKA Grant No. K68322, and the CNRS while
this author visited Laboratoire G-SCOP, Grenoble}\\
\small Computer and Automation Research Institute\\[-0.8ex]
\small Hungarian Academy of Sciences\\[-0.8ex]
\small Budapest, P.O. Box 63, H-1518, Hungary\\[-0.8ex]
\small \texttt{gyarfas@sztaki.hu}\\[-0.8ex]\\
\and Andr\'as Seb\H o\\
\small CNRS, Laboratoire G-SCOP, \\[-0.8ex]
\small 46 avenue F\'elix Viallet, 38031 Grenoble Cedex, France \\[-0.8ex]
\small \texttt{andras.sebo@g-scop.inpg.fr}\\ \\
 \and Nicolas
Trotignon\\
\small CNRS, LIAFA, Universit\'e Paris 7, Paris Diderot\\[-0.8ex]
 \small 175 rue du Chevaleret, 75013 Paris, France\\[-0.8ex]
\small \texttt{nicolas.trotignon@liafa.jussieu.fr}}

\newcommand{\sm}{\setminus}
\newcommand{\gap}{{\rm{gap}}}

\begin{document}
\maketitle
\newpage
\begin{abstract}

The {\em chromatic gap} is the difference between the chromatic
number and the clique number of a graph.  Here we investigate
$\gap(n)$, the maximum chromatic gap over graphs on $n$ vertices.
Can the extremal graphs be explored? While computational problems
related to the chromatic gap are hopeless, an interplay between
Ramsey theory and matching theory leads to a simple and (almost)
exact formula for $\gap(n)$ in terms of Ramsey numbers. For our
purposes it is more convenient to work with the {\em covering
gap}, the difference between the clique cover number and stability
number of a graph and this is what we call the {\em gap} of a
graph. Then $\gap(n)$ can be equivalently defined (by switching
from a graph to its complement), as the maximum gap over graphs of
$n$ vertices. Notice that the well-studied family of perfect
graphs are the graphs whose induced subgraphs have gap zero.
Our study is a first step towards better understanding of graphs whose
induced subgraphs have gap at most $t$. The maximum of the (covering) gap and
the chromatic gap running on all induced subgraphs
will be called {\em perfectness gap}.

Using $\alpha(G)$ for the cardinality of a largest stable
(independent) set of a graph $G$, we define $\alpha(n)=\min
\alpha(G)$ where the minimum is taken over triangle-free graphs on
$n$ vertices. It is easy to observe that $\alpha(n)$ is
essentially an inverse Ramsey function, defined by the relation
$R(3,\alpha(n))\le n < R(3,\alpha(n)+1)$. Our main result is that
$\gap(n)=\lceil n/2\rceil - \alpha(n)$, possibly with the exception
of small intervals (of length at most $15$) around the Ramsey
numbers $R(3,m)$, where the error is at most 3.

The central notions in our investigations are the {\em gap-critical}
and the {\em gap-extremal} graphs. A graph $G$ is
\emph{gap-critical} if for every proper  induced subgraph $H\subset G$,
$\gap(H) < \gap(G)$ and {\em gap-extremal} if it is gap-critical
with as few vertices as possible (among gap-critical graphs of the
same gap). The strong perfect graph theorem, solving a long
standing conjecture of Berge that stimulated a broad area of
research, states that gap-critical graphs with gap $1$ are the holes
(chordless odd cycles of length at least five) and antiholes
(complements of holes). The next step, the complete description of
gap-critical graphs with gap $2$ would probably be a very difficult
task. As a
very first step, we prove that there is a unique $2$-extremal graph,
$2C_5$, the union of two disjoint (chordless) cycles of length five.

In general, for $t\geq 0$, we denote by $s(t)$ the smallest order of a
graph with gap $t$ and we call a graph is \emph{$t$-extremal} if it
has gap $t$ and order $s(t)$. It is tempting to conjecture that
$s(t)=5t$ with equality for the graph $tC_5$. However, for $t\ge 3$
the graph $tC_5$ has gap $t$ but it is not gap-extremal (although
gap-critical). We shall prove that $s(3)=13$, $s(4)=17$ and $s(5)\in
\{20,21\}$. Somewhat surprisingly, after the uncertain values $s(6)\in
\{23,24,25\},$ $s(7)\in \{26,27,28\}$, $s(8)\in\{29,30,31\}$, $s(9)\in
\{32,33\}$ we can show that $s(10)=35$.  On the other hand we can
easily show that $s(t)$ is asymptotically equal to $2t$, that is,
$\gap(n)$ is asymptotic to $n/2$. According to our main result the gap
is actually equal to $\lceil n/2 \rceil - \alpha(n)$, unless $n$ is
in an interval $[R, R+14]$ where $R$ is a Ramsey number, and if this
exception occurs the gap may be larger than this value by only a small
constant (at most $3$).

The definition of $s(t)$ does not change if we replace (covering) gap by chromatic gap,
so it can in fact be defined with the perfectness gap as well: it is the smallest order
of a graph with perfectness gap equal to $t$.

Our study provides some new properties of Ramsey-graphs themselves: it shows that triangle-free Ramsey graphs have high
matchability and connectivity properties, therewith providing some new
properties of the Ramsey-graphs themselves, and leading possibly to
new bounds on Ramsey-numbers.
\end{abstract}

\section{Introduction}

After the proof of the strong perfect graph conjecture \cite{CST}, the
problems concerning graph families that are close to perfectness
become more interesting.  Here we focus our attention on a parameter
that we call  the {\em chromatic gap} of a graph,
$\gap(G)$,
equal to the ``duality
gap'' of a most natural integer linear programming formulation of the
graph coloring problem.

Graphs in this paper are  undirected, their vertex set is denoted by
$V(G)$. A {\em cycle} is a connected subgraph with all degrees equal
to $2$.  A {\em clique} is a subset of the vertices inducing a
complete subgraph, and a  {\em stable set} does not induce any edge.
The notations $C_i$ and $K_i$ will refer to cycles, respectively cliques of order $i$ $(i=1,2,\ldots)$ .

The size of a largest clique (resp.\ stable set) in a graph $G$ is
denoted, by $\omega(G)$ (resp.\ $\alpha(G)$).  We also speak about
$k$-cliques or $k$-stable sets meaning that their cardinality is
$k$. A $3$-clique is also called a {\em triangle}.  The \it chromatic
number, \rm $\chi(G)$, and \it clique-cover number, \rm $\theta(G)$,
denote the minimum number of partition classes of $V(G)$ into stable
sets and into complete subgraphs, respectively.  Using $\overline G$
for the complement of $G$, we have obviously

\begin{equation}\label{complement}
\omega(G)=\alpha(\overline G), \chi(G)=\theta(\overline G)
\end{equation}
and
\begin{equation}\label{basicineq}
\chi(G)\ge \omega(G)\ge {|V(G)|\over \theta(G)},  \theta(G)\ge
\alpha(G)\ge {|V(G)|\over \chi(G)}.
\end{equation}

Let us define the \it chromatic gap \rm of a graph $G$ as
$\chi(G)-\omega(G)$, and the \it covering gap \rm as
$\theta(G)-\alpha(G)$.  Although these parameters are equivalent
(through (\ref{complement})), for our purposes it is more
convenient to work with the latter, so we define the {\em gap, or
covering gap} of a graph $G$ \rm as
$\gap(G)=\theta(G)-\alpha(G)$. Notice that {\em perfect graphs}
are the graphs whose induced subgraphs have gap zero. The {\em perfectness gap}
of a graph is the maximum of the (covering) gap and
the chromatic gap running on all induced subgraphs.

A graph $G$ is \emph{gap-critical} if for every proper   induced
subgraph $H\subset G$, $\gap(H) < \gap(G)$. The perfect graph theorem
\cite{CST} states that gap-critical graphs with gap $1$ are the holes
(chordless odd cycles of length at least five) and antiholes
(complements of holes). The complete description of gap-critical
graphs with gap $2$ would probably be a very difficult task - it seems
there is not even a plausible guess available. Trivial members can be
obtained as a disjoint union of holes and/or antiholes. A nontrivial
member(15 vertices, $\alpha=6, \theta=8$) is shown in \cite{PRSURR},
p. 427. Deleting any pair of vertices of the Ramsey-graph $R_{13}$,
the unique graph with $\omega(G)=2, \alpha(G)=4$ on $13$ vertices,
gives another example of order $11$ with $\alpha=4,
\theta=6$. However, as we shall prove, the smallest order of a
gap-critical graph with gap $2$ is 10, the unique example is the
trivial member, the union of two disjoint $C_5$. The graph $R_{13}$
itself is also gap-critical with gap $3$, in fact the smallest one
(see Section \ref{sub:start}).

Note that the definition of gap-critical graphs cannot be simplified
by requiring only $\gap(G-v)<\gap(G)$ for every vertex $v$: indeed,
for instance the gap of the circular graph $C(3,3)$ (on $10$
cyclically ordered vertices, where any three cyclically consecutive
ones form a clique) is $1$, deleting any vertex the gap is $0$
although $C_5$ subgraphs are still present.  Here the smallest
example.  Consider a hole on 5 vertices $c_1 \dots c_5 c_1$ and
replicate $c_1$ and $c_3$ (replicating a vertex $v$ means adding a
vertex adjacent to $v$ and all neighbors of $v$).  For the obtained
graph $G$ we have $\omega(G)=3$, $\chi(G)=4$, but for any $v\in V(G)$
$\omega(G-v) = \chi(G-v) = 3$, while $G$ contains a $C_5$. So, the
complementary graph $\bar G$ is not gap-critical, although $\gap (G-v)
< \gap (G)$ for all $v\in V(G)$.

The central topic of our work is to determine the maximum gap of
graphs of order $n$, denoted by $\gap(n)$ which leads to a study of
gap-extremal graphs. For $t\geq 0$, we denote by $s(t)$ the smallest
order of a graph with gap $t$.  A graph is \emph{$t$-extremal} if it
has gap $t$ and order $s(t)$; it is {\em gap-extremal}, if it is
$t$-extremal for some $t$. Note that the empty graph has gap 0, so
$s(0) = 0$, and --- since $C_5$ is the unique smallest non-perfect
graph --- $s(1) = 5$, and $C_5$ is the only $1$-extremal graph. It
will be much more difficult to prove that $s(2)=10$
(Theorem~\ref{thm:s210}). It is tempting to conjecture that the
pattern continues and $s(t)=5t$ with equality for the graph $tC_5$,
this is how we started $\ldots$ However, classical Ramsey-graphs
provide better bounds. We shall prove that $s(3)=13, s(4)=17$ and
$s(5)=21$ or $20$. From a general conjecture we think that the true
value is $21$. Somewhat surprisingly, after the uncertain values
$s(6)\in \{23,24,25\}, s(7)\in \{26,27,28\}, s(8)\in\{29,30,31\},
s(9)\in \{32,33\}$ we can show that $s(10)=35$.

Gap-extremal graphs are obviously gap-critical. Holes and antiholes
are gap-critical but if they have more than five vertices they are
not gap-extremal; if they have more than eleven vertices their gap
is also not maximal among graphs of the same order, since the gap of
two disjoint $C_5$ is $2$.

\medskip
A large $\theta(G)$ might be the consequence of a small $\omega(G)$. But small
clique number may mean not too many edges, so a large $\alpha(G)$
too! What happens with the gap in this competition? The trade
between the size of cliques and stable sets is described by
Ramsey-theory, itself having a lot of open questions. We will
convert the relations provided by Ramsey numbers into a balance
between $\theta$ and $\alpha$. Using Ramsey-numbers as a black box
we will be able to (almost) determine our  functions.

It will turn out to be essentially true that the graphs with a large
gap are triangle free. In other words, decreasing the clique-size,
makes $\theta$ increase more than it does increase $\alpha$. To work
out this precisely will need a refined analysis based on details
concerning Ramsey-numbers $R(3,.)$ and matchings. In
Section~\ref{sec:mind} we prove simple statements about the gap, about
Ramsey-numbers and about matchings that will provide the right tools
for this work. In Section~\ref{sec:find} we determine the gap function
with only a small constant error, and this relies mainly on a study of
triangle-free graphs.

In view of this role, we will need to use variants of the notions
and terms for triangle-free graphs separately. We will speak about
{\em triangle-free $t$-extremal graphs} which means that their
cardinality is minimum among triangle-free graphs of gap $t$. Note
that a triangle-free gap-extremal graph is not necessarily a
gap-extremal graph, since there might be a graph containing a
triangle with smaller cardinality and the same gap. By analogy, the
corresponding notations for triangle-free graphs will be
$\gap_2(n)$, $s_2(t)$. Thus $\gap_2(n)$ is the maximum gap among
triangle-free graphs on $n$ vertices, $s_2(t)$ is the smallest order
of a triangle-free graph with gap $t$. Clearly, $\gap(n)\ge
\gap_2(n)$ for all $n\in\Nset$, and $s(t)\le s_2(t)$ for all
$t\in\Nset$. $(\Nset$ is the set of natural numbers $\{1,2,\ldots\})$.

\medskip For any $n\in\Nset$, $t=\gap(n)$, adding $n-s(t)$ isolated
points to a $t$-extremal graph we get a graph of maximum gap among
graphs of order $n$. However, both $\alpha$ and $\theta$ increase by
the addition of isolated vertices. When $G$ is triangle free, graphs
of maximum gap, at the same time with minimum stability number among
triangle-free graphs on $n$ vertices will be particularly appreciated.
Let $\alpha(n)$ denote the minimum of $\alpha(G)$ over triangle-free
graphs $G$ with $n$ vertices.  So, $\alpha(n)$ is defined by the
relation $R(3,\alpha(n))\le n < R(3,\alpha(n)+1)$.  A graph $G$ on $n$
vertices will be called {\em stable gap-optimal}, if $G$ is triangle
free, $\gap(G)=\gap_2(n)$, and $\alpha(G)=\alpha(n)$. It will turn out
that there exist stable gap-optimal graphs for every $n$.  Therefore
it is unavoidable to know something about the function $\alpha(n)$, in
fact it is just the inverse of the well studied Ramsey function
$R(3,x)$.

We say that a graph  is an  {\em $(\omega,\alpha)$-Ramsey graph} $
(\omega,\alpha \in\Nset)$ if it is of maximum order among the graphs
$G$ without an $\omega$-clique (a clique of size $\omega$) and
without an $\alpha$-stable set (stable set of size $\alpha$). By
Ramsey's theorem \cite{LEX}, this maximum is finite. The   smallest
$n$ such that for any graph $G$ of order $n$ either
$\omega(G)\ge\omega$ or $\alpha(G)\ge \alpha$, is called the  Ramsey
number  $R(\omega,\alpha)$. We will use mainly Ramsey numbers for
$\omega=3$. Clearly, the order of $(\omega,\alpha)$-Ramsey graphs is
$R(\alpha,\omega)-1,$ and their maximum clique and stable set have
size $\omega-1$, $\alpha -1$.

Clearly, the above introduced number $\alpha(n)$ $(n\in\Nset)$ is
actually defined by the relation $R(3,\alpha)\le n < R(3,\alpha
+1)$. It is equal to {\em the number of Ramsey-numbers smaller than or
  equal to $n$.}  Indeed, among the Ramsey-numbers $R(3,x)$ those with
$1,2,\ldots x$ are smaller than or equal to $n$, and all the others
are larger. It will turn out that $s(t+1) - s(t)$ is usually $2$, and
the exceptions are at the Ramsey-numbers where this difference is
equal to $4$ with rare exceptions 5 of 3 (but these latter might
actually all be for $t\le 3$).

Although $s(t)$ will be determined with a constant error (modulo
Ramsey numbers), we also include a transparent easy proof in Section
2 that shows that $2t+c_1\sqrt{tlog{t}}\le s(t)\le
2t+c_2\sqrt{tlog{t}}$ (Corollary \ref{gapass}).

The main result of the paper is finding $\gap(n)$ and $s(t)$ with
constant error in terms of Ramsey numbers. First we shall prove that
$\gap(n)=\gap_2(n)=\lceil n/2\rceil - \alpha(n)$ except when $n$ is
even and there exists odd numbers $n_1,n_2$ such that $n=n_1+n_2$ and
$\alpha(n)=\alpha(n_1)+\alpha(n_2)$, in which case $1$ must be
added. The exceptional case can occur in an obvious way, when $n$ is
a Ramsey number and $n_1$ or $n_2$ is equal to 1, or in a rather
mysterious way (only if $n_1=n_2=5 ?)$, when we call $n$ Ramsey-perfect.

A number $n$ is \emph{Ramsey-perfect} if $n$ is not an even
Ramsey-number and $n=n_1 + n_2$, where $n_1, n_2 \geq 5$ are odd and
$\alpha(n) = \alpha(n_1) +\alpha(n_2)$.   We know only one
Ramsey-perfect number, 10 ($\alpha(10)=2\alpha(5)$), and we believe
that there are no others. One way this might still happen is
$\alpha(n)=\alpha(n-5)+\alpha(5)$, in that case $n-1,n-4$ must be
both Ramsey numbers --- we call them {\em (Ramsey) twins}. Probably
there are no Ramsey twins beyond $6,9$ but this is not proved,
although Erd\H os and S\'os \cite{ESO} (see also in \cite{CG})
conjectured $R(3,m+1)-R(3,m)$ tends to infinity with $m$. Our main
results are summarized as follows.

\begin{itemize}

\item[--]   $\gap_2(n)=\lceil n/2\rceil - \alpha(n) +\varepsilon(n)$,
where $\varepsilon(n)=1$ if $n$ is an even Ramsey-number or a
Ramsey-perfect number and $0$ otherwise (Theorem \ref{thm:gap2}).

\item[--] The functions $\gap(n), s(t)$ are
determined with a small error by their restricted counterparts:  for
all $n,t\in\Nset$: $0\le \gap(n) - \gap_2(n) \le 2,\quad$ $0\le
s_2(t) - s(t) \le 10.$ (Theorem \ref{thm:bounds}).

\item[--] A synthesis of this work: for all $n\in {\Nset}\setminus \cup
_{\alpha\in \Nset}[R(3,\alpha), R(3,\alpha)+14]:$
$\gap(n)=\gap_2(n)=\lceil n/2\rceil - \alpha(n)$, and always $\lceil
n/2\rceil - \alpha(n)\le \gap(n) \le   \lceil n/2\rceil - \alpha(n)
+ 3.$ (Theorem \ref{thm:main}).
\end{itemize}

It is worth noting that for Ramsey numbers $R$ that are at least $5$
bigger than the preceding Ramsey-number (so maybe for all
Ramsey-numbers larger than $28$), only one $s(t)$ value is uncertain
and equal to either $R+1$ or $R+2$.  Also, our study reveals high
matchability and connectivity properties of Ramsey graphs. For
example, $(3,\alpha+1)$-Ramsey-graphs are $(R(3,\alpha+1)-R(3,\alpha)
-3)$-connected, moreover, deleting at most $R(3,\alpha+1)-R(3,\alpha)
-3$ vertices, the remaining $n\ge R(3,\alpha) +2$ vertices, if $n$ is
even, induce a graph with a perfect matching (Corollary
\ref{cor:matching}).

Finally we mention some related works. B\'\i r\'o \cite{biro} raised
the related problem of finding the minimum of $\alpha$ while fixing
$n$ and $\theta$, more precisely finding
$$\beta(n,\theta)=\min\{\alpha (G) :  G \hbox{ graph, } |V(G)|=n,
\theta(G)=\theta\}$$ and gave the first bounds and a conjecture.
Jahanbekam and West \cite{west} stated another conjecture for
constrained values of $n$ and $\theta$. If
$\theta \ge {n+1\over 2}$ Theorem~\ref{thm:gap2} easily provides the following formula for $\theta$, implying these conjectures:
 $\beta(n,\theta)=n+\alpha(W)-W - \varepsilon$, where $W=2(n-\theta)+1$ and $\varepsilon$ is $0$ or $1$, the latter if $W$ is Ramsey-perfect or another (even more exceptional, possibly non-existing) case that we will neglect here.   A recent communication of
B\'\i r\'o, F\"uredi and Jahanbekam, \cite{BFJ} proves a formula for
$\beta(n,\theta)$ in the range $\theta \ge {n+3\over 2}$ with similar methods\footnote{\cite{BFJ}
mentions the relation of $\beta(n,\theta)$ to the present work (to
\cite{leibniz} or to an
 earlier version from November 2009), notice an inaccuracy, but miss the close tie to
 Theorem~\ref{thm:gap2}.}. The equality between the two formulas can be proved easily (the different appearance of \cite{BFJ} is due to exploiting less inequalities between Ramsey numbers).  As far as we know, finding the exact values of $\gap(n)$ (without restricting ourselves to triangle-free graphs) and
the solution of B\'\i r\'o's problem {\em for arbitrary $\theta$ }
both remain open problems.

\section{Asymptotic of $s(t)$}

Before giving the exact values of the function $\gap$ and $\gap_2$
(up to a small constant) we show how to get easily the asymptotic
of $s(t)$.

\begin{proposition}\label{upperbound}\pont
$s(t)\le s_2(t)\le 2t+ c_1\sqrt{t\log{t}}$.
\end{proposition}

\noindent \bf Proof. \rm The celebrated result of Kim  \cite{KIM}
states that for every sufficiently large $n$ there is a graph $G_n$
with $n$ vertices such that $\omega(G_n)=2$ and $\alpha(G_n)\le
9\sqrt{n\log{n}}.$ Define $f(t)$ as the smallest $n$ for which there
exists $G_n$ such that
\begin{equation}\label{RAM}
\left\lceil{n\over 2}\right\rceil -9\sqrt{n\log{n}}\ge t.
\end{equation}
Clearly $f(t)$ is an upper bound for $s_2(t)$ because by the
definition of $G_n$ and by (\ref{RAM})
\begin{equation}\label{RAM1}
\gap(G_n)=\theta(G_n)- \alpha(G_n)\ge \left\lceil {n\over
2}\right\rceil - 9\sqrt{n\log{n}}\ge t.
\end{equation}
One can easily check that the last inequality in (\ref{RAM1}) can be
satisfied with $n=2t+\lfloor c_1\sqrt{t\log{t}}\rfloor$ where $c_1$ is
a constant. This gives the required upper bound.
 \qed

\begin{proposition}\label{ebed}\pont $s(t)\ge 2t + \alpha(2t)\ge 2t+c_2\sqrt{t\log{t}}.$
\end{proposition}

\noindent \bf Proof. \rm Let $G$ be a graph with $\gap(G)=t$ and with
$n$ vertices. Consider a clique cover of $G$ obtained by greedily
selecting a largest clique in the subgraph induced by the vertex set
uncovered in previous steps. Suppose that in the first $k$ steps
cliques of size at least three were selected, covering $2k$ vertices
plus a set $A\subseteq V(G)$, followed by $l$ steps of selecting
edges and covering $Y$, finally a set $Z$ of independent vertices
covers the rest of the vertices of $G$. Set $B=Y\cup Z$.

Then clearly,
$$\theta(G)\le {n-|Z|-|A| \over 2} +|Z| = {n-|A| \over 2} +{|Z| \over 2}\le {n-|A|\over 2}+{\alpha(B)\over 2}$$
therefore
$$\theta(G)-\alpha(G)\le {n-|A|\over 2}+{\alpha(B)\over 2}-\alpha(G)\le
{n-|A|\over 2}+{\alpha(B)\over 2}-\alpha(B) = {n-|A|-\alpha(B)\over
  2}$$

thus $2t + |A|+ \alpha(B)\le n=s(t)$.  We gained $|A|+\alpha(B)$
over the 2t lower bound. However, we know that $3|A|+|B|\ge n \ge
2t$. It is easy to see that the gain is smallest for $|A|=0$ thus we
gain at least $\alpha(2t)$ as desired. \qed

\begin{corollary}\label{gapass} $s(t)=2t+\theta(\sqrt{t\log{t}})$
\end{corollary}

\section{Matchings and Ramsey numbers}\label{sec:mind}

In this section we explore the main properties of the gap of a
graph, of gap-critical graphs, of the relation of these to matchings
and the Ramsey-numbers.

\subsection{Easy facts}\label{subsec:easy}

\begin{proposition}
  \label{l:components}
  If a graph $G$ has $k$ connected components $C_1, \dots, C_k$ then
  $\gap(G) = \gap (C_1) + \cdots + \gap(C_k)$.  Every connected
  component of a gap-critical graph is gap-critical. Every connected
  component of a gap-extremal graph is gap-extremal.
\end{proposition}

\noindent
\bf Proof\rm : Both $\theta$ and $\alpha$ are sums of the $\theta$ and $\alpha$ of the components. \qed

\medskip\noindent

\begin{proposition}
  \label{p:gapmonotonicity} The $\Nset\rightarrow \Nset$ functions
  $\gap$ and $\gap_2$ are monotone increasing.
\end{proposition}

\noindent \bf Proof. \rm Indeed, if $n_1\le n_2$, then adding
$n_2-n_1$ isolated vertices to a graph $G$ of order $n_1$ of maximum
gap, we get a graph of order $n_2$ of the same gap. \qed

\begin{proposition}
  \label{p:gapsubadditivity} For any $n_1,n_2\in\Nset$ we have
  $\gap(n_1+n_2)\ge \gap(n_1)+gap(n_2).$ For any $t_1,t_2\in\Nset$ we
  have $s(t_1+t_2)\le s(t_1)+s(t_2).$
\end{proposition}

\medskip\noindent
\bf Proof \rm: Let $G$ be a graph that consists of two components,
$G_1$  on $n_1$ vertices, and $G_2$ on $n_2$ vertices,
$\gap(G_1)=\gap(n_1)$ and $\gap(G_2)=\gap(n_2)$. Then $G$ has
$n_1+n_2$ vertices, and $\gap(n_1+n_2)\ge
\gap(G)=\gap(n_1)+\gap(n_2)$. For the second part let $G$ be a
graph that consists of two components,  a $t_1$-extremal graph $G_1$
on $s(t_1)$ vertices, and  a $t_2$-extremal graph $G_2$ on $s(t_2)$
vertices. Then $G$ has $s(t_1)+s(t_2)$ vertices, and
$\gap(G)=t_1+t_2$ thus $s(t_1+t_2)\le |V(G)|=s(t_1)+s(t_2)$. \qed

\medskip The equality is easily satisfied, for instance $\gap(5)=1$,
$\gap(17)=4$, and $\gap(22)=5$ as we will see in
Section~\ref{sub:start} . We have a third, similar inequality where
the condition of equality is less trivial (Theorem~\ref{thm:perfect}),
that turns out to be very restrictive and the related notion of
Ramsey-perfect numbers are crucial for the main results
(Subsection~\ref{sub:s2}).
\begin{proposition}\label{p:subadditivity} For any $n_1,n_2\in\Nset$ we have $$\alpha(n_1+n_2)\le \alpha(n_1)+\alpha(n_2).$$
 \end{proposition}

\medskip\noindent
\bf Proof \rm: Indeed, a graph $G$ that consists of two components,
$G_1$  on $n_1$ vertices,
 and $G_2$ on $n_2$ vertices, $\alpha(G_1)=\alpha(n_1)$ and $\alpha(G_2)=\alpha(n_2)$,
 has $n_1+n_2$ vertices, and $\alpha(n_1+n_2)\le
 \alpha(G)=\alpha(n_1)+\alpha(n_2)$.\qed

\begin{proposition}
  \label{l:removeK}
  Let $G$ be a graph and $Q$  a clique of
  $G$.  Then
  $$\theta(G) \ge \theta(G - Q) \ge \theta(G) - 1,$$

  $$\alpha(G) \ge \alpha(G - Q) \ge
  \alpha(G)-1, \,$$

  $$\gap(G)+ 1\ge \gap(G - Q) \ge  \gap(G)- 1,$$

  and there exists a chain of induced subgraphs of $G$ with gaps equal to
  $\gap(G), \gap(G) -1,\ldots, 0.$ Furthermore, if $G$ is
  gap-critical,
  $$\theta(G - Q) = \theta(G) - 1,\, \alpha(G - Q) =
  \alpha(G),\,\, \gap(G - Q) = \gap(G)- 1.$$
\end{proposition}

Notice that the equality $\gap(G - Q)=\gap(G)-1$ may hold also for
graphs that are not gap-critical (see the example in the Introduction:
a hole on 5 vertices with two non-adjacent vertices replicated).

\medskip\noindent \bf Proof\rm: $\theta(G) \le \theta(G - Q) + 1$ is
true because adding $Q$ to any clique cover of $G-Q$ we get a clique
cover of $G$.  $\alpha(G)\le \alpha(G - Q) + 1$ holds because any
stable set meets $Q$ in at most one vertex.  The third inequality
follows from these first two and the obvious bounds
$\alpha(G-Q)\le\alpha(G)$, $\theta(G-Q)\le\theta(G)$. The statement
about the chain of induced subgraphs follows by noting that the deletion of a
vertex changes the gap by at most $1$, in the beginning it is
$\gap(G)$, and at the end it is $0$.

If $G$ is gap-critical,
 $\gap(G-Q)=\gap(G) + 1$,  $\gap(G-Q)=\gap(G)$
 cannot occur in the proven inequalities,
  so the only option is $\gap(G- Q)
  = \gap(G) - 1$, and then $\theta(G - Q) =
  \theta(G) - 1$ and $\alpha(G - Q) = \alpha (G)$.
\qed

\medskip
A vertex of a graph is \emph{simplicial  } if its neighbors
induce a complete graph.

\begin{proposition}
  \label{simplicial}
  If $G$ is  gap-critical, then  it has no simplicial vertex.
\end{proposition}

\noindent
\bf Proof. \rm
 If $v \in V(G)$ is a simplicial vertex,
 $\alpha(G - N[v]) = \alpha (G) -1,$ since $S\cup \{v\}$ is a stable-set for any
 stable set $S$ of $G - N[v]$, contradicting
 Proposition~\ref{l:removeK} for $Q=N[v]$.
\qed

\medskip
The following generalizes the condition on $N(v)$ if $\alpha\le 2$:
\begin{proposition}\label{l:alpha2}
Let $G$ be a graph such that $\alpha(G)\le 2$ and there exists $v\in V(G)$ where $G(N(v))$ is perfect.
Then $\gap (G)\le 1.$
\end{proposition}

\noindent \bf Proof\rm: Consider $G_1:=G(N[v])$ which is now perfect,
and $Q:=G- N[v]$ which is a clique because of $\alpha(G)\le 2$.  By
Proposition~\ref{l:removeK}, $0 = \gap(G_1) = \gap(G-Q) \ge \gap(G) -
1$.  \qed

At last we state easy but crucial lower bounds for $s(t)$ and $s_2(t)$, and an interesting relation between
these bounds and the equality $s(t)=s_2(t)$.

\begin{proposition}
  \label{l:jumpk}
  If there exists a $(t+1)$-extremal graph $G$ with $\omega(G)\ge k$ $(k\in\Nset)$, then
  $s(t+1) \geq s(t) +k,$ in particular, for any $t\in\Nset$:
   $s(t+1) \geq s(t) +2$.
\end{proposition}

\noindent
\bf Proof\rm:
  Let $K$ be a $k$-clique in $G$. By Proposition~\ref{l:removeK}, $\gap(G\sm K) = \gap(G) - 1$.  So $$s(t)
  \leq |V(G\sm K)| = |V(G)| - k = s(t+1) - k.\hbox{\,\,\,\qquad\qquad\qed}$$

\medskip
We prove three simple but important statements on the relation of $s$ and $s_2$ :

\begin{proposition}
  \label{sess2} If $s(t+1)=s(t)+2$, then $s(t)=s_2(t)$, $s(t+1)=s_2(t+1)$.

  If  $s(t)\ne s_2(t)$ or $s(t+1)\ne s_2(t+1),$ then $\,\,\, s(t+1)\ge\,\, s(t)+3$.

  If $s(t)=s_2(t)$ and $\quad s(t+1)\ne s_2(t+1),$ then   $s_2(t+1)\ge s_2(t)+4$.
\end{proposition}

\medskip\noindent
\bf Proof\rm: Let $G$ be $t+1$-extremal, and suppose $s(t+1)=s(t)+2$. If $G$ is not triangle-free,
by Proposition~\ref{l:removeK},
 $s(t+1)\ge s(t)+3$,
so $G$ is triangle-free, and deleting the two endpoints of an edge, the gap decreases by $1$,
so what we get is $t$-extremal,
and the first statement follows. The second statement is just the indirect reformulation of the first.
The third follows by
$s_2(t+1)>s(t+1)\ge s(t)+3=s_2(t)+3,$ using the preceding inequality.
\qed

\subsection{Gaps and Matchings}\label{subsec:matchings}

As usual,  $\nu(G)$  denotes the size of a maximum matching of $G$,
the maximum number of pairwise disjoint edges; let $\zeta(G)$
denote the minimum number of edges that cover the vertices of $G$.
If $G$ is a triangle-free graph, $\theta(G)=\zeta(G)$.
 The reader can find in any textbook or check
that for connected graphs $\nu(G) +\zeta(G)=n.$

A graph is \emph{factor-critical} if the removal of any vertex yields
a graph with a perfect matching. (It is convenient to include graphs
of order $1$ under this term.)  A graph is {\em bicritical} if
deleting any two vertices there is a perfect matching. Clearly, {\em
  factor-critical and bicritical graphs are connected}.  The following
is a simple but ingenious and important result of Gallai
\cite{gallai:factorCritical} (in English in \cite{LP} or \cite{LEX}
Exercise 26 page 58).

\begin{theorem}[Gallai, \cite{gallai:factorCritical}]
  \label{th:gallai}
  If $G$ is connected and  $\nu(G\sm v) = \nu(G)$ for all $v \in
  V(G)$, then $G$ is factor-critical, and in particular it has an odd number of vertices.
\end{theorem}

\begin{proposition}
  \label{l:fcrit}
  If $G$ is a triangle-free and gap-critical graph then every component
  of $G$ is factor-critical of (odd) order at least $5$.
  \end{proposition}

\noindent
\bf Proof\rm: Let $H$ be a component of a triangle-free, gap-critical graph.  By
  Proposition~\ref{l:components} $H$ is gap-critical.
  Since $H$ is triangle-free, $\theta(H) = \zeta(H)$ and by
  Proposition~\ref{l:removeK}, for all $v\in V(H)$ we have $\zeta (H\sm v)=\theta (H\sm v)
  =\theta(H) - 1 = \zeta(H) - 1$. So
  $$\nu(H\sm v) = |V(H\sm v)| - \zeta(H\sm v) = |V(H)| - \zeta(H) =
  \nu(H),$$
whence $H$ is factor-critical by Theorem~\ref{th:gallai}.

If some component is a vertex, deleting that isolated vertex the
gap does not decrease. It cannot be a triangle either. \qed

\medskip
The following proposition gives a lower bound on the gap and this
bound will turn out to be {\em very sharp}, in fact an equality. The
intuition behind it: {\em in a triangle-free graph $G$
$\theta(G)=\theta(G-v)+1$ for every vertex $v\in V(G)$} implies
$\theta(G) = \lceil {V(G)\over 2} \rceil,$ which is the smallest
possible value in a triangle-free graph. That is, if we want
$\theta(G)$ to be largest possible comparing to $\theta(G-v)$, then
$\theta$ takes its smallest possible value.

\begin{proposition}\label{p:lower} For any triangle-free graph $G$,
$\gap(G)\ge \lceil {|V(G)| \over 2 }\rceil- \alpha(G),$ and for
connected triangle-free gap-critical graphs the equality holds. If
there exists a triangle-free gap-extremal graph of order $n$ with
$k$ components of order $n_1,\ldots,n_k$,
$$ \gap_2(n)= \lceil {n_1 \over 2 }\rceil -\alpha(n_1) +\ldots +\lceil {n_k \over 2 }\rceil -\alpha(n_k).$$
\end{proposition}

\noindent \bf Proof\rm:  Since $G$ is triangle-free, $\theta(G)\ge
\lceil {|V(G)| \over 2 }\rceil$ so $\gap(G)=\theta(G)-\alpha(G)\ge
\lceil {|V(G)| \over 2 }\rceil- \alpha(G)$. If $G$ is gap-critical
and connected, by Proposition~\ref{l:fcrit} it is factor-critical,
so $\theta(G)=\lceil {|V(G)|\over 2}\rceil$, settling the first
claim.
Now if $G$ is triangle-free gap-extremal, then by
Proposition~\ref{l:components} all of its components are connected
gap-critical graphs, and by the already proven assertion,
$\gap(G_i)=   \lceil {n_i \over 2 }\rceil -\alpha(G_i).$

If $\alpha(G_i) > \alpha(n_i)$ then replacing $G_i$ by $H_i$ of the
same order $n_i$, triangle-free, ($\theta(H_i)\ge \lceil
n_i\rceil$), and $\alpha(H_i)=\alpha(n_i)<\alpha(G_i)$, the gap
increases,
contradicting that $G_i$ is gap-extremal. So $\theta(G_i)=\lceil
{n_i\over 2}\rceil$, $\alpha(G_i)=\alpha(n_i)$, finishing the
proof with an application of Proposition~\ref{l:components}. \qed

\medskip

Is the triangle-free condition essential in these statements ? For
some of the claims it can be dropped! Gallai himself proved in
\cite{gallai:colorCritical}: {\em If the complement of a
$k$-color-critical graph is connected, it has at least $2k-1$
vertices.} By Proposition~\ref{l:removeK} the complements of
gap-critical graphs are color-critical, so we immediately get:

\begin{proposition}\label{p:gallaicol} If $G$ is a connected gap-critical graph,
$\theta(G) \le \lceil { |V(G)|\over 2 } \rceil.$
\end{proposition}

Stehl\'\i k \cite{ST} proved the sharpening of Gallai's general
theorem stating that {\em there exists a coloration where all color
  classes are of size at least two}, extending Gallai's proof
\cite{LEX}, \cite{LP} of Theorem~\ref{th:gallai}
\cite{gallai:factorCritical}. Despite these promising generalizations,
we were not able to make essential use of
Proposition~\ref{p:gallaicol} or prove in any other way that
gap-extremal graphs cannot contain a triangle. However,
Proposition~\ref{sess2}, the main results of the paper and further
verifications for small $t$ (see Section~\ref{sub:start}) suggest that
it is true:

\begin{conjecture}\label{conj:trianglefree} Every gap-extremal graph is triangle-free.
\end{conjecture}

\subsection{Gaps and Ramsey-numbers}\label{subsec:ramsey}

Let $W_8$ be the Wagner' graph \cite{RA}, a cycle on $8$ vertices with
its four long chords. Deleting one of these chords we get $W_{81}$ and
deleting two neighboring chords we get $W_{82}$.  Let $R_{13}$ be the
graph on $\{r_1, \dots, r_{13}\}$ with the following edges:
$r_ir_{i+1}$ and $r_ir_{i+5}$, $i=1, \dots, 13$, where the addition is
taken modulo 13. It is well known \cite{RA} that $R_{13}$ is the
largest graph such that $\omega=2$ and $\alpha=4$.  Note that
$\gap(R_{13}) = 3$.

The following is mostly an extract  of \cite{RA}, except for the
lower bounds on $R(3,24),$ $\ldots, R(3,29)$ that are from
\cite{WSLX}:

\begin{proposition}\label{p:firstramsey} The Ramsey-numbers $R(3,l)$ for values
$l=2, 3, 4, 5, 6, 7, 8, 9$ are $3$, $6,$ $9,$ $14,$ $18,$ $23,$
$28,$ $36$, and the corresponding Ramsey graphs are unique for
$l=2$, $l=3$ and $l=5$: $K_2$, $C_5$  and $R_{13}$ respectively. For
$l=4$ there are three Ramsey-graphs, $W_8$, $W_{81}$, $W_{82}$.
Moreover $40\le R(3,10)\le 43$, $46\le R(3,11)\le 51$, $52\le
R(3,12)\le 59$, $59\le R(3,13)\le 69$, $66\le R(3,14)\le 78$, $73\le
R(3,15)\le 88$, and  $R(3,16)\ge 79$,  $R(3,17)\ge 92$, $R(3,18)\ge
99$,  $R(3,19)\ge 106$,   $R(3,20)\ge 111$, $R(3,21)\ge 122$,
$R(3,22)\ge 125$, $R(3,23)\ge 136$,   $R(3,24)\ge 143$, $R(3,25)\ge
153$, $R(3,26)\ge 159$, $R(3,27)\ge 167$, $R(3,28)\ge 172$,
$R(3,29)\ge 182.$

\smallskip\noindent
$R(4,4)=18$, and the unique $(4,4)$-Ramsey graph on $17$ vertices is
a cycle of length $17$ with all chords between vertices at distance
$2$, $4$, $8$.
 \end {proposition}

\medskip

The following is a result of Xiaodong, Zheng and Radziszowski
\cite{XXR} (Theorem 3) see also \cite{RA} 2.3 (g).

\begin{proposition}\label{thm:conn}\cite{XXR}
 If $p, q \ge 2$, $R(3,p+q-1)\ge R(3,p) + R(3,q) + \min\{p,q\} - 2.$
\end{proposition}

 \begin{proposition}\label{p:iconsecutive}
 \begin{itemize}
 \item[ \petit{cons1}]
 $\alpha+1 \ge R(3,\alpha + 1) - R(3,\alpha)\ge 3$ (provided $\alpha\ge 2$ for the second inequality) and both
 inequalities are strict if
 both $R(3,\alpha)$ and $R(3,\alpha+1)$ are even.
  \item[\petit{cons2}] $R(3,\alpha +2 ) - R(3,\alpha)\ge 7$ provided $\alpha\ge 3$.
  \item[\petit{cons3}] $R(3,\alpha +3 ) - R(3,\alpha)\ge 11$ provided $\alpha\ge 2$.
  \item[\petit{cons4}] $R(3,\alpha +4 ) - R(3,\alpha)\ge 17$ provided $\alpha\ge 3$.
  \item[\petit{cons5}] $R(3,\alpha +k ) - R(3,\alpha)\ge R(3,k+1) +k -1$, if $\alpha\ge k+1\ge 3$.
   \item[\petit{cons55}] The right hand side of (5) for $\alpha\ge 3$ and $k=5,6,7$   are: $22, 28, 34.$
    \item[\petit{cons56}] The right hand side of (5) for $\alpha\ge 4$, $k=8, 9, 10, 11, 12, 13$
    are: $43, 48, 55, 62, 70, 78.$
  \item[\petit{cons6}] $R(3,\alpha +14) - R(3,\alpha)\ge 86$, if $\alpha \ge 3$.
    \end{itemize}
\end{proposition}

\medskip\noindent \bf Proof\rm: First, we prove (\ref{cons1}): The
upper bound is the easy and most well-known upper bound $R(3,\alpha+1)
\le R(3,\alpha) + R(2,\alpha+1)$ \cite{LEX}, where the equality does
not hold if both terms on the right hand side are even, and where of
course $R(2,\alpha+1)=\alpha+1$ (to see this, start the usual
induction with a vertex of even degree).  Since equality would imply
that $R(2,\alpha+1)=\alpha+1$ is even too (that is, $\alpha$ is odd),
we have the assertion concerning the upper bound.  The lower bound of
(\ref{cons1}) is a result in \cite{BEFS} and also a special case of
Proposition~\ref{thm:conn} by substituting $q=2$ and $R(3,2)=3$.

Second, we check (\ref{cons2}) by substituting $p=\alpha\ge 3$, $q=3$
and $R(3,3)=6$ into Proposition~\ref{thm:conn}. Third, substituting
$p=\alpha\ge 4$, $q=4$ and $R(3,4)=9$ into Proposition~\ref{thm:conn}
provides (\ref{cons3}) for $\alpha\ge 4$, and for $\alpha=2,3$ it can
be checked in Proposition~\ref{p:firstramsey}.  (\ref{cons4}) for
$\alpha\ge 5$ is a specialization, and can be checked directly in
Table~I for $\alpha=3,4$, (\ref{cons5}) is just a rewriting of
Proposition~\ref{thm:conn}.

Finally, if we specialize (\ref{cons5}) to $k=5,..., 14$, we get
(\ref{cons55}), (\ref{cons56}), (\ref{cons6}) for $\alpha\ge 6,
\ldots, \alpha\ge 15$, respectively. For $\alpha=3,\ldots,9$ we still
get the inequalities from \cite{RA} Table II and I, for the lower
bounds are provided until $l=23$, and the upper bounds until $l=15:$
for instance, $R(3,23)\ge 136$, $R,(3,9)=36$, so $R(3,23) - R,(3,9)\ge
100.$ For the lower bounds concerning the highest arguments we have to
rely on upper bounds \cite{WSLX} copied into
Proposition~\ref{p:iconsecutive}.  The inequalities with the largest
values that we have to check are $R(3,\alpha+14) - R(3,\alpha)\ge 86,$
for $\alpha=4,\ldots, 14$. (For $\alpha\ge 15$ we have from
(\ref{cons5}) and substituting $R(3,15)\ge 73$ from Table I \cite{RA}
$R(3,\alpha+14) - R(3,\alpha)\ge R(3,15)+ 13\ge 86$.)  We make the
last checking, for $\alpha=14$: $R(3,28) - R(3,14)\ge 86.$ Indeed,
from Proposition~\ref{p:firstramsey} $R(3,28)\ge 172$ (copied from
\cite{WSLX}) and $R(3,14)\le 78$ (from
Proposition~\ref{p:firstramsey}), so in fact $R(3,28) - R(3,14)\ge
94\ge 86.$ \qed

\medskip If $R(3,\alpha+1) - R(3,\alpha)=3$, we will say that
$R(3,\alpha)$, $R(3,\alpha+1)$ are {\em twins}.

\begin{proposition}\label{p:lowern}
$\gap_2(n)\ge \lceil {n \over 2 }\rceil- \alpha(n).$
\end{proposition}

\bf Proof\rm: Indeed, by Proposition~\ref{p:lower}  for any
triangle-free graph $G$ on $n$ vertices $\gap_2(n)\ge \gap(G)\ge
\lceil {n \over 2 }\rceil- \alpha(G),$ and if we apply this to a
triangle-free graph $G$ with $\alpha(G)=\alpha(n)$ we get the
claim. \qed

\bigskip
We will now need  to deduce conditions on the  equality in
Proposition~\ref{p:subadditivity}. These computations will enable us
to conclude that there exist  stable gap-optimal graphs of any order
$n\in\Nset$, and this will be crucial for our formulas describing
the gap. A combination of the inequalities of
Proposition~\ref{thm:conn} and the upper bound of
Proposition~\ref{p:iconsecutive} (\ref{cons1}) yield the following
characterization of the equality in
Proposition~\ref{p:subadditivity} that will be crucial for
describing the gap-function, through Ramsey-perfect numbers.

\begin{theorem}\label{thm:perfect} Let $n, n_1, n_2, n_3\in\Nset$. Equality in $\alpha(n_1+n_2)\le \alpha(n_1)+\alpha(n_2)$ implies
  that there exist $\varepsilon, \varepsilon_1, \varepsilon_2$ such
  that $n_1+n_2 - \varepsilon$, $n_1 + 1 + \varepsilon_1$,
  $n_2+1+\varepsilon_2$ are all Ramsey-numbers, and
  $\varepsilon,\varepsilon_1, \varepsilon_2\in\{0,1\}$, $\varepsilon +
  \varepsilon_1 + \varepsilon_2\le 1.$

  Furthermore if  $n_i\ge 3$ for $i=1,2,3$ then
  $\alpha(n_1+n_2+n_3) < \alpha(n_1)+\alpha(n_2) + \alpha(n_3).$

\end{theorem}

In the last, strict inequality the condition is necessary:
$\alpha(6)=3=3\alpha(2);$ if say \marginpar{!} $n_3=2$, then
$n:=n_1+n_2+n_3$ may be a Ramsey number, $n-3$ its twin, and $n-2$
could be Ramsey-perfect. However, luckily, we are interested in
these equalities only if the numbers $n_1$, $n_2$, $n_3$ are odd,
and then a stronger inequality holds (Lemma~\ref{lem:ideal}).

Note that even in the first part of the theorem,
 $\alpha(n_1+n_2)=\alpha(n_1)+\alpha(n_2)$ with $n_2=1$ can be useful. This holds if and only if $n_1+1$ is a Ramsey-number.
 If in addition $n_1+1$ is even, a Ramsey-graph on $n_1$ vertices and an isolated vertex provides the maximum gap
  (Theorem~\ref{thm:gap2}).


\medskip
\noindent \bf Proof\rm : We  reprove the easy  inequality
$\alpha(n_1+n_2)\le \alpha(n_1)+\alpha(n_2)$ (see
Proposition~\ref{p:subadditivity}) in a complicated way, in order to
deduce the conditions of equality. Set $\alpha_i=\alpha(n_i)$. Then
$n_i\le R(3,\alpha_i+1)-1$ $(i=1,2)$.
\begin{lemma}\label{lem:perfect}
For arbitrary  $\alpha_1, \alpha_2\in\Nset$
$$R(3,\alpha_1+1)-1+ R(3,\alpha_2+1)-1\le R(3,\alpha_1+\alpha_2) + 1, \hbox{ and}\leqno\petit{lem3}$$
equality implies that Proposition~\ref{p:iconsecutive} (\ref{cons1})
(first part) holds with equality for the smaller of $\alpha_1,
\alpha_2$.
\end{lemma}

\noindent \bf Proof\rm : By symmetry we may suppose
$\alpha_1\ge\alpha_2$.

If $\alpha_2=1$ then (\ref{lem3}) and
Proposition~\ref{p:iconsecutive} (\ref{cons1}) (first part) are
equalities. If $\alpha_2\ge 2$ we can substitute $p=\alpha_1+1$,
$q=\alpha_2$ into Proposition~\ref{thm:conn} and add $1$ to both
sides :
$$R(3,\alpha_1+1)-1+ R(3,\alpha_2)-1 + \alpha_2 +1\le R(3,\alpha_1+\alpha_2)+1. \leqno\petit{lem1}$$
Applying Proposition~\ref{p:iconsecutive} (\ref{cons1}) to
$\alpha_2$,
$$R(3,\alpha_2) + \alpha_2 +1 \ge R(3,\alpha_2+1), \leqno\petit{lem2}$$
and (\ref{lem1}),(\ref{lem2})gives lemma (together with the remark
on equality). \qed

From the definitions and from Lemma \ref{lem:perfect}, $n_1+n_2\le
R(3,\alpha_1+1)-1+ R(3,\alpha_2+1)-1\le R(3,\alpha_1+\alpha_2) + 1,$
from where we indeed can read $\alpha(n_1+n_2)\le
\alpha_1+\alpha_2$, and the equality holds if and only if
$$R(3,\alpha_1+\alpha_2)  \le n_1+n_2 \le R(3,\alpha_1+\alpha_2) + 1 .$$
These inequalities allow at most one of $n_1$ or $n_2$ be one less
than $R(3,\alpha_1+1)-1$ or $R(3,\alpha_2+1)-1$ respectively, that
is, $\varepsilon_1+\varepsilon_2\le 1$,  and in case of equality,
$n_1+n_2=R(3,\alpha_1+\alpha_2)$, that is, $\varepsilon=0$.

\medskip

Next we prove the second part of Theorem \ref{thm:perfect}, the
strict inequality when $n$ is decomposed into three numbers. We
could apply Lemma \ref{lem:perfect} twice and each time the
conditions for the equality in it, but then the result we get would
be too weak. We repeat the proof, applying
Proposition~\ref{thm:conn} directly, twice, choosing its arguments
carefully:

\begin{lemma}\label{lem:three} For arbitrary natural numbers $\alpha_1\ge\alpha_2\ge\alpha_3\ge 2$,
$$R(3,\alpha_1+1)-1+ R(3,\alpha_2+1)-1+R(3,\alpha_3+1)-1\le R(3,\alpha_1+\alpha_2+\alpha_3-1) +2.\leqno{\petit{threexact}}$$
\end{lemma}

\medskip
Lemma~\ref{lem:three} concludes the proof of
Theorem~\ref{thm:perfect} since $n_1+n_2+n_3$ is less than or equal
to the left hand side of (\ref{threexact}). Since  $n_i\ge 3$
implies $\alpha_i\ge 2$, Lemma~\ref{lem:three} shows that
$n_1+n_2+n_3$ is also bounded from above by the right hand side of
(\ref{threexact}). Then, because of Proposition~\ref{p:iconsecutive}
(1) (second inequality providing the lower bound $3$), the right
hand side can be upper bounded by $R(3,\alpha_1+\alpha_2+\alpha_3) -
1,$ proving that $\alpha (n_1+n_2+n_3) \le
\alpha_1+\alpha_2+\alpha_3 -1,$ showing the claimed strict
inequality of Theorem~\ref{thm:perfect}. \qed

\medskip\noindent
{\bf Proof of Lemma~\ref{lem:three}}:
Apply the upper bound of Proposition~\ref{p:iconsecutive}  (1) to get that the left hand side is less than or equal to
$$R(3,\alpha_1+1)-1+(R(3,\alpha_2 )+ \alpha_2 + R(3,\alpha_3 )+ \alpha_3),\leqno{\petit{icons}}$$
where the sum in the parentheses can in  turn be bounded according to Proposition~\ref{thm:conn}:
$$R(3,\alpha_2 )+ R(3,\alpha_3 )+ \alpha_2 + \alpha_3\le R(3,\alpha_2+\alpha_3 -1)- (\alpha_3-2) +\alpha_2 +\alpha_3.\leqno{\petit{seize}}$$
Substituting this to (\ref{icons}) and applying Proposition~\ref{thm:conn} again to the result, (\ref{icons})$\le$
$$R(3,\alpha_1+1) +R(3,\alpha_2+\alpha_3 -1) +\alpha_2 +1\le R(3,\alpha_1+\alpha_2+\alpha_3-1)-
(\alpha_2+1-2) +\alpha_2 +1,$$ after noting that
$\alpha_2+1\le\min\{\alpha_1+1,\alpha_2+\alpha_3-1\}.$ \qed


\smallskip
\begin{corollary}\label{cor:perfeq} A number $n\in\Nset$ is Ramsey-perfect if and only if
there exist $\alpha_1\ge\alpha_2\ge 2$  that satisfy
$n=R(3,\alpha_1+\alpha_2)+1= R(3,\alpha_1+1)-1 + R(3,\alpha_2+1)-1$, where $R(3,\alpha_i+1)$
is even $(i=1,2)$. Moreover, then the equality holds in (\ref{lem1}), (\ref{lem2}).
\end{corollary}

\noindent \bf Proof\rm: Indeed, if $n$ is Ramsey-perfect, let
$n_1,n_2\ge 5$ be odd numbers such that $n=n_1+n_2$,
$\alpha(n)=\alpha(n_1)+\alpha(n_2)$. Since $n_1$ and $n_2$ satisfy the
condition of Theorem~\ref{thm:perfect}, the theorem can be
applied. Denote $\alpha:=\alpha(n)$, $\alpha_1:=\alpha(n_1)\ge 2$,
$\alpha_2:=\alpha(n_2)\ge 2$. Since $n$ is not a Ramsey-number,
$\varepsilon=1$, and then $\varepsilon_1=\varepsilon_2=0.$ In other
words $n=R(3,\alpha)+1$, $n_i=R(3,\alpha_i+1)-1$ are odd, $(i=1,2),$
$n=n_1+n_2$, $\alpha=\alpha_1+\alpha_2,$ showing the assertion.
Moreover, Lemma \ref{lem:perfect} is satisfied with equality, whence
(\ref{lem1}), (\ref{lem2}) as well. Conversely, if the equality and
the parity condition are satisfied with $\alpha_1, \alpha_2\ge 2$,
then defining, $n_1:=R(3,\alpha_1+1)-1$, $n_2:=R(3,\alpha_1+1)-1$ we
see that $n=R(3,\alpha_1+\alpha_2)$ is Ramsey-perfect. \qed

\bigskip
The lack of other examples of twins or other Ramsey-perfect numbers is not really surprising: only the first nine
Ramsey values are known. Yet we believe that all the applied inequalities cannot be
tight for arbitrary large Ramsey-numbers, so
we state two conjectures:
\begin{conjecture}   The natural number $n$ is Ramsey-perfect if and only if $n$ is even and $n-1$ is
the bigger of  Ramsey twins.
\end{conjecture}

\begin{conjecture} The only  Ramsey twins are $\{3,6\}$ and $\{6,9\}$.
\end{conjecture}

\begin{corollary}\label{cor:extremal}
  Let $G$ be triangle-free-extremal with a a minimum number of
  components.  Then $G$ has at most two components, and two if and
  only if $n:=|V(G)|$ is Ramsey-perfect, when $$\gap(G)= \lceil
  n/2\rceil -\alpha(n)+1,$$ otherwise $n$ is odd, $G$ is connected,
  and $$\gap(G)= \lceil n/2\rceil -\alpha(n).$$ In both cases the
  triangle-free-extremal graphs are stable gap-optimal, and in the
  second case any triangle-free graph on $n$ vertices and stability
  number $\alpha(n)$ is stable gap-optimal.
\end{corollary}

\noindent
\bf Proof\rm:  Let $G$ be a triangle-free-$t$-extremal graph with a minimum number of components,
$t\in\Nset$, and let $G_1,\ldots,G_k$ be its components, of order $n_1,\ldots, n_k$,  $n:=|V(G)|=n_1+\ldots, n_k$.
According to Proposition~\ref{l:fcrit} all the components are factor-critical, in particular all the $n_i$ are odd,
$\theta(G_i)=\lceil n_i/2\rceil$ $(i=1,\ldots, k)$, and by Proposition~\ref{p:lower},
$$\gap(G)=\lceil n_1/2\rceil - \alpha(n_1)+ ... + \lceil n_k/2\rceil - \alpha(n_k).\leqno{\petit{gap}}$$

It follows now from Theorem~\ref{thm:perfect} that $k\le 2$, because otherwise three components can be replaced by one,
contradicting the choice of $G$:
$$\left\lceil {n_1+n_2+n_3\over 2}\right\rceil\ge \lceil n_1/2\rceil + \lceil n_2/2\rceil + \lceil n_3/2\rceil-1,\,\,
 \alpha(n_1+n_2+n_3)\le \alpha(n_1)+\alpha(n_2) + \alpha(n_3) -1.$$
Two components can also be replaced by just one, unless the equality is satisfied in both of the following inequalities:
 $$\left\lceil {n_1+n_2\over 2}\right\rceil\ge  \lceil n_1/2\rceil + \lceil n_2/2\rceil-1,\,\,
 \alpha(n_1+n_2)\le \alpha(n_1)+\alpha(n_2).$$
So $k=1$, or $k=2$, and then (\ref{gap}) specializes to the claimed formula, since for $k=2$
$$\gap(G)= \lceil n_1/2\rceil -\alpha(n_1) + \lceil n_2/2\rceil-\alpha(n_2)=\left\lceil {n_1 +n_2\over 2}
\right\rceil+1 - \alpha(n_1+n_2),$$ and this happens if and only if
$n$ is Ramsey-perfect.

In both cases $G$ is stable gap-optimal, and conversely, if
$n=s_2(t)$ is neither an even Ramsey-number nor Ramsey-perfect, then
according to  Proposition~\ref{p:lower} every graph $H$ on $n$
vertices and stability number $\alpha(n)$ satisfies: $\gap(G)\ge
\lceil n/2\rceil - \alpha(n)=\gap_2(n),$ so there is equality
throughout, and $G$ is stable gap-optimal. \qed

\section{Finding the gap with constant error}\label{sec:find}

In this section we first determine the functions $\gap_2(n)$ and
$s_2(t)$ exactly, and then the functions $\gap(n)$ and $s(t)$ with
small errors ($2$ and $10$ respectively), moreover we prove that the
error may occur only after Ramsey-numbers on an interval of length
$13$.

\subsection{Finding the triangle-free gap}\label{sub:s2}

Recall that $\gap_2(n)$ is the maximum of the gap of a triangle-free
graph of order $n$, and $s_2(t)$ denotes the minimum order of a
triangle-free graph of gap $t$. The main result of this section is a
simple formula for these functions if the inverse Ramsey numbers
$\alpha(n)$ are used as black boxes.

\begin{theorem}\label{thm:gap2} $\gap_2(n)=\lceil n/2\rceil - \alpha(n) +\varepsilon(n)$,
where $\varepsilon(n)=1$ if $n$ is an even Ramsey-number, or if it is Ramsey-perfect, and $0$ otherwise.
\end{theorem}



\medskip
\noindent
{\bf Proof}: Let $f(n):=\lceil n/2\rceil - \alpha(n)+\varepsilon(n)$.

\medskip\noindent
{\bf Claim~1}: $\gap_2(n)\ge f(n)$  for all $n\in\Nset.$

\smallskip Indeed, if $n$ is neither an even Ramsey-number nor Ramsey-perfect, this is just
Proposition~\ref{p:lowern}. If $n$ is an even Ramsey-number, then $\alpha(n-1)=\alpha(n)-1$ and
$\lceil {n-1\over 2}\rceil = \lceil {n\over 2}\rceil$, so by the monotonicity of $\gap_2$
(see Proposition~\ref{p:gapmonotonicity}):
$$\gap_2(n)\ge \gap_2(n-1)\ge f(n-1)= \lceil n/2\rceil - \alpha(n)+1.$$
More generally, if $n=n_1+n_2$ where $n_1$, $n_2$ are odd numbers  and $\alpha(n)=\alpha(n_1)+\alpha(n_2)$, then
$\lceil {n\over 2}\rceil +1= \lceil {n_1\over 2}\rceil+ \lceil {n_2 \over 2}\rceil$,
and applying  Proposition~\ref{p:gapsubadditivity} and then Proposition~\ref{p:lowern}:
$$\gap_2(n)\ge \gap_2(n_1)+ \gap_2(n_2)\ge \lceil {n_1\over 2}\rceil - \alpha(n_1) +
\lceil {n_2\over 2}\rceil - \alpha(n_2)= \lceil n/2\rceil - \alpha(n)+1.$$

Corollary~\ref{cor:extremal} establishes the theorem for the values
$n=s_2(t)$  $(t=1, 2,\ldots )$, thus we get

\smallskip\noindent
{\bf Claim~2}: If $n=s_2(t)$ for some $t\in\Nset$, then $\gap_2(n)=f(n).$

\medskip\noindent
{\bf Claim~3}: The function $f(n)$ is monotone increasing.

\medskip Indeed, since $\lceil n/2\rceil$ is a monotone increasing function, we have $f(n+1)\ge f(n)$ unless
$\alpha(n)$ is increasing, or unless $\varepsilon(n)$ is decreasing when $n$ grows to $n+1$.
We prove that in both of these less trivial events actually $f(n+1)=f(n)$:

Assume first that $\alpha(n+1)> \alpha(n)$. Then $\alpha(n+1)= \alpha(n)+1$, that is, $n+1$ is the
Ramsey-number $R(3,\alpha(n)+1)$. If in addition $n$ is even,
$\lceil {n+1\over 2}\rceil=\lceil {n\over 2}\rceil+1$, and $\varepsilon(n+1)= 0=\varepsilon(n)$
since $n+1$ is an odd Ramsey-number, so neither $n$ nor $n+1$ is an even Ramsey-number or
Ramsey-perfect by Theorem~\ref{thm:perfect}. So
$$f(n+1) = \lceil {n+1\over 2}\rceil - \alpha(n+1)+\varepsilon(n+1)= \lceil {n\over 2}\rceil +1 -
(\alpha(n)+1)+\varepsilon(n)+0= f(n).$$
If $n$ is odd -- and still $\alpha(n+1)> \alpha(n)$ -- , then $\lceil {n+1\over 2}\rceil=\lceil {n\over 2}\rceil$,
but then $n+1$ is an even Ramsey-number, so
$$f(n+1) = \lceil {n+1\over 2}\rceil - \alpha(n+1)+\varepsilon(n+1)= \lceil {n\over 2}\rceil  -
(\alpha(n)+1)+\varepsilon(n)+1= f(n).$$

Second, assume that $\alpha(n+1)= \alpha(n)$, but
$\varepsilon(n+1)=\varepsilon(n)-1.$ Then $\varepsilon(n)=1$, so $n$
is even, and therefore $\lceil {n+1\over 2}\rceil=\lceil {n\over
2}\rceil+1$, so again $f(n+1)=f(n)$ proving the claim.

\medskip To finish the proof of the theorem,
suppose for a contradiction that $\gap_2\ne f$. Let $x$ be the
smallest integer $x$ for which $t:=\gap_2(x)\ne f(x)$. By Claim~1,
$\gap_2(x)> f(x).$ Then, by Claim~3, we have for all $y\le x$:  $
t=\gap_2(x)>f(x)\ge f(y)=\gap_2(y)$ by the minimality of $x$.   So
$s_2(t)=x$, and then, by Claim~2, $\gap_2(x)=f(x)$, a contradiction
that proves the theorem. \qed

\begin{corollary}\label{cor:Ramseygap} For all $\alpha\in\Nset,$ $\gap_2(R(3,\alpha) ) = \lceil
{R(3,\alpha)+1\over 2}\rceil - \alpha=\gap_2(R(3,\alpha) -1),$ in particular, Ramsey-numbers are not
in the image of the function $s_2$.
\end{corollary}

\medskip
\noindent
{\bf Proof}: If $n$ is even, $\varepsilon(n)=1$, so  $\lceil n/2\rceil - \alpha +\varepsilon(n)= \lceil {n+1\over 2}\rceil -
\alpha$. If $n$ is odd, $\varepsilon(n)=0$ and $\lceil {n\over 2}\rceil=\lceil {n+1\over 2}\rceil$, so again
$\lceil n/2\rceil - \alpha +\varepsilon(n)= \lceil {n+1\over 2}\rceil - \alpha$.
In both cases $\gap_2(n)=\gap_2(n-1)$, so $n\ne s_2(t)$ for any $t$.
\qed

\begin{corollary}For every $\alpha\in\Nset$ for which $R(3,\alpha+1)-R(3,\alpha)\ge 4$, exactly the odd
numbers of the interval $[R(3,\alpha) +3, R(3,\alpha+1)-1]$ are the values of the function $s_2(t)$,  for
$t=\gap_2(R(3,\alpha) )+1,\ldots, \gap_2(R(3,\alpha+1))-1$.
\end{corollary}

\noindent
{\bf Proof}: This is an immediate consequence of Theorem~\ref{thm:gap2}, since for the integers $n$ of the given
interval both $\alpha(n)$ and $\varepsilon(n)$ are constant, and $\lceil n/2\rceil$ increases exactly on odd numbers.
\qed

\begin{corollary}\label{cor:alphagapopt} For every $n$ there exists a  stable gap-optimal graph $G$, defined from
an arbitrary $(3,\alpha+1)$-Ramsey graph  $G_\alpha$ $(\alpha=1,2, \ldots):$
\begin{itemize}
\item[--] if $n\in [R(3,\alpha)+2, R(3,\alpha+1)-1]$  or if $n=R(3,\alpha)+1$ is not Ramsey-perfect or if
$n=R(3,\alpha)$ is odd, let $G$ be an arbitrary,  order $n$ induced subgraph of $G_\alpha$.

\item[--] if $n$ is Ramsey-perfect, $n=R(3,\alpha)+1=n_1+n_2$,
$n_i:=R(3,\alpha_i+1)-1$ is odd $(i=1,2)$,
$\alpha=\alpha_1+\alpha_2$,  then let $G$ consist of two
components:  $G_{\alpha_1}$ and $G_{\alpha_2}$. \item[--] if
$n=R(3,\alpha)$ is even, let $G$ consist of $G_{\alpha -1}$ and an
isolated vertex.
\end{itemize}
If $n$ or $n-1$ is equal to $R(3,\alpha)$  then $G$ is not
necessarily connected, but otherwise every stable gap-optimal graph
is connected.
\end{corollary}

For $n=6$ the only stable gap-optimal graph is $C_5$ and an
isolated vertex. For $n=7$ and any number $R(3,\alpha)+1$ which is
not Ramsey-perfect, a graph having two components, a Ramsey-graph
and a $K_2$ is stable gap-optimal, and may actually coincide with
$G_\alpha$.

\medskip
\noindent {\bf Proof}: In the first case $\gap(G)=\lceil n/2\rceil
- \alpha(n)=\gap_2(n)$ according to Theorem \ref{thm:gap2} $G$ is
indeed stable gap-optimal.

In the second and third case, if $n=R(3,\alpha)+1$ or $n=
R(3,\alpha)$, the defined graphs are readily stable gap optimal,
and so are the graphs of the remark before the proof if
$n=R(3,\alpha)+1$ but $n$ is not Ramsey-perfect. If $n$ is neither
of these two numbers, it cannot be written as the sum of two
nonzero numbers whose inverse Ramsey-numbers sum up to $\alpha(n)$
(see Theorem~\ref{thm:perfect}), so the defined stable gap-optimal
$G$ is connected. \qed

\begin{corollary}\label{cor:matching} $(3,\alpha+1)$-Ramsey-graphs are
$(R(3,\alpha+1)-R(3,\alpha) -3)$-connected, moreover, deleting at most $R(3,\alpha+1)-R(3,\alpha) -3$
vertices, the remaining $n\ge R(3,\alpha) +2$ vertices, if $n$ is even, induce a graph with a perfect matching.
\end{corollary}

\medskip
\noindent {\bf Proof}: Apply Corollary~\ref{cor:alphagapopt} to
odd $n\in [R(3,\alpha) +3, R(3,\alpha+1)-1]$ : any induced subgraph $G$ of
$G_\alpha$ on $n$ vertices has optimal gap. Fix this $G$ and have
a look at Theorem~\ref{thm:perfect}:
$\varepsilon(n)=\varepsilon(n-1)=0$, and we see that the
jump-points of the function $\gap_2(n)=\lceil n/2\rceil - \alpha$,
that is the values of the function $s_2(t)$ on the considered
interval are exactly the odd numbers. So $G$ is a
triangle-free-extremal graph, and either by
Corollary~\ref{cor:alphagapopt} or by Corollary~\ref{cor:extremal}
it is connected, and by Proposition~\ref{l:fcrit} it is
factor-critical, and the graphs in the assertion arise by deleting
a vertex in such a graph.\qed

\begin{corollary}
Order $n$ induced subgraphs of $(3,\alpha+1)$-Ramsey-graphs induce a factor-critical graph if
$n\ge R(3,\alpha) +3$ is odd, and a bicritical graph if $n\ge R(3,\alpha) +4$ is even. \qed
\end{corollary}

Indeed, this corollary is an immediate consequence of Corollary~\ref{cor:matching}.

\medskip We now determine the recurrence relations for the function $s_2$. Why?
Doesn't Theorem~\ref{thm:gap2} tell us all we need? Indeed, it does
already tell the most important information, the following theorem
and its proof are secondary, the reader can skip it at first
reading.  However, besides an automatic conversion of
Theorem~\ref{thm:gap2} from the $\gap_2$ function to $s_2$, it also
has a new content: it shows that for a Ramsey-perfect number $n$,
the interval $[n,n+3]$ cannot contain a Ramsey-number again. Besides
making the formulas simpler (at the price of a slightly more
difficult proof), it reveals some interesting relations between the
distance of consecutive Ramsey numbers and Ramsey-perfectness.

\begin{corollary}\label{thm:s2}  For all $t$, $s_2(t)$ is odd or
  Ramsey-perfect. Moreover, the function $s_2$ is determined by the
  following recursive relations:
 \begin{itemize}
 \item[1.]  If  neither $s_2(t)+1$, nor $s_2(t)+2$ are Ramsey, then

                 $1.1\,\, s_2(t+1)=s_2(t)+2$ if $s_2(t)$ is not Ramsey-perfect.

                 $1.2\,\, s_2(t+1)=s_2(t)+3$ if $s_2(t)$ is Ramsey-perfect, and $s_2(t)+3$ is not Ramsey.

                 $1.3\,\, s_2(t+1)=s_2(t)+4$ if $s_2(t)$ is Ramsey-perfect, $s_2(t)+3$ is Ramsey,

                 \qquad \quad\qquad\qquad\qquad\qquad moreover $s_2(t)+4$ is Ramsey-perfect.

                 $1.4\,\, s_2(t+1)=s_2(t)+5$ if $s_2(t)$ is Ramsey-perfect, $s_2(t)+3$ is Ramsey,

                 \qquad \quad\qquad\qquad\qquad\qquad but $s_2(t)+4$ is not Ramsey-perfect.

                 \item[2.]  If   either $s_2(t)+1$  or $s_2(t)+2$ are Ramsey, then

                 $2.1\,\,  s_2(t+1)=s_2(t)+3,$ if $s_2(t)+3$ is Ramsey-perfect.

                 $2.2\,\, s_2(t+1)=s_2(t)+4$ otherwise, except if $s_2(t)+4$ is   Ramsey.

                 $2.3\,\, s_2(t+1)=s_2(t)+5,$ if $s_2(t)+4$ is  Ramsey.
 \end{itemize}
\end{corollary}



\medskip
\noindent {\bf Proof}: Let $n=s_2(t)$ then by definition,
$\gap_2(n)>\gap_2(n-1)$, and let $\alpha:=\alpha(n),
\varepsilon=\varepsilon(n).$ Suppose that $n$ is odd, or
Ramsey-perfect. We will show that the recursive relations 1.1-2.3
hold, and $s_2(t+1)$ is also odd or Ramsey-perfect.

\medskip\noindent
1.1:  If  neither $n+1$, nor $n+2$ are Ramsey-numbers, and $n$ is
not Ramsey-perfect, then by assumption $n$ is odd and  $\alpha$,
$\varepsilon$ are constant in the interval $[n,n+2]$. Therefore by
Theorem~\ref{thm:gap2} $\lceil {n\over 2}\rceil=\lceil {n+1\over
2}\rceil<\lceil {n+2\over 2}\rceil$, so 1.1 holds.

\medskip\noindent
1.2: If $n=s_2(t)$ is Ramsey-perfect then according to
Corollary~\ref{cor:perfeq} there exist $\alpha, \alpha_1,
\alpha_2\in\Nset$ such that
$$n=R(3,\alpha)+1=R(3,\alpha_1+1)-1 + R(3,\alpha_2+1)-1, \hbox{ $n$ is even.}\leqno{\petit{eq:perfect}}$$
According to Theorem~\ref{thm:gap2},
$\gap_2(n)=\gap_2(n+1)=\gap_2(n+2)$, since while the ceiling
increases by $1$, $\varepsilon$ decreases by $1$. Now
$\gap_2(n+3)=\gap_2(n)+1$ unless $n+3$ is a Ramsey-number again, and
1.2 is checked.

\medskip\noindent
1.3: If $n+3$ is a Ramsey-number (and otherwise the same condition holds as in 1.2), then in addition
to $\gap_2(n)=\gap_2(n+1)=\gap_2(n+2)$ we have $\gap_2(n+2)=\gap_2(n+3)$, since both $\theta$ and $\alpha$ have increased.
However,  $n+4$ may or may not be Ramsey-perfect, and in the former case  $\gap_2(n+4)=\gap_2(n)+1$, that is,
$s(t+1)=n+1$, as claimed.

\medskip\noindent
1.4:  In case $n+4$ is not Ramsey-perfect (and otherwise the same condition holds as in 1.2) $\gap_2(n+3)=\gap_2(n+4)$ and
$n+4$ is even, so $\theta$, $\alpha$ remain the same as for $n+3$. However, $n+5$ is odd, and cannot be Ramsey again
since $n+3$ is Ramsey; $\theta$ increases, but $\alpha$ does not:  $\gap_2(n+5) > \gap_2(n)= \gap_2(n+4)$, so $s_2(t+1)=n+5$,
as claimed.

\medskip\noindent
2.1: If $n+3$ is Ramsey-perfect then $n+2$ is an odd
Ramsey-number, and by Theorem~\ref{thm:gap2} we have by parity,
and because of
$\varepsilon(n)=\varepsilon(n+1)=\varepsilon(n+2)=0$,
$\varepsilon(n+3)=1$, $\alpha(n)=\alpha(n+1)=\alpha$,
$\alpha(n+2)=\alpha(n+3)=\alpha+1$:
$\gap_2(n)=\gap_2(n+1)=\gap_2(n+2)<\gap_2(n+3)$ as claimed.

\medskip\noindent
2.2: If the same hold but $n+3$ is not Ramsey-perfect, then all
the relations of 2.1 hold except that we have now
$\varepsilon(n+3)=0$, and therefore
$\gap_2(n)=\gap_2(n+1)=\gap_2(n+2)=\gap_2(n+3)< \gap_2(n+4),$
where $n+4$ is indeed odd.

\medskip\noindent
2.3: If $n+1$ or $n+2$ is a Ramsey-number, and $n+4$ is a
Ramsey-number again, then $n$ is odd, $n+1$ and $n+4$ are
twins. So $n+1$ is an even Ramsey-number, and
$\alpha(n+1)=\alpha(n)+1$, $\varepsilon(n+1)=\varepsilon(n)+1=1$
compensate one another, so Theorem~\ref{thm:gap2} gives this time
$\gap_2(n)=\gap_2(n+1)=\gap_2(n+2)=\gap_2(n+3)=\gap_2(n+4)<\gap_2(n+5).$
Note that $s_2(t+1)=n+5$ is even in this case, in accordance with
the fact that $n+5$ is Ramsey-perfect because of
$\alpha(n+5)=\alpha(n)+\alpha(5)=\alpha +2$.
\qed


Corollary \ref{thm:s2} gives concrete values of $s_2(i)$ for $i<12$,
because we do not know whether 40 or 41 is a Ramsey number.

\begin{corollary}\label{values} The values of $s_2(i),i=1,\dots,11$ are $5,10,13,17,21,25,29,31,33,35,39.$
\end{corollary}

In fact, we will prove $s_2(i)=s(i)$ almost everywhere, and we
conjecture it is true everywhere. This is a slightly weaker
conjecture than Conjecture~\ref{conj:trianglefree}.

\begin{conjecture}\label{conj:gapgap2}
$\gap(n)= \gap_2(n)$ for all $n\in\Nset$, and $s(t)= s_2(t)$  for all $t\in\Nset$.
\end{conjecture}

In the next section we show that the possible exceptions to this
conjecture are  at constant distance from Ramsey-numbers, and at any
such place the difference of the function value from the ``usual''
$\lceil n/2\rceil -\alpha(n)$ is also a small constant.

\subsection{Bounding the gap function}\label{sub:s}


The first assertion of the following lemma states that once the
relation $s(t)=s_2(t)$ holds, it  surely holds again and again
(together with the equivalent equality $\gap(t)=\gap_2(t)$) until
the next Ramsey-number; the second assertion ensures that the
relation $s(t)=s_2(t)$ holds again after exceptions restricted to a
small interval (of size at most $29$) after each Ramsey-number.

\begin{lemma}\label{thm:constant} Assume $R(3,\alpha)\le s(t)=s_2(t) < R(3,\alpha+1).$  Then
\begin{itemize}
\item[--]For all $t'\in\Nset$ such that
$s(t) \le s(t') \le R(3,\alpha+1): \,\, s(t')=s_2(t')$
\item[--]There exists $t'\in\Nset$, $t  <  t'\le t+29$ such that
$$s(t) = n < R(3,\alpha+1)<s(t')=s_2(t')\le R(3,\alpha+1) + 85\le R(3,\alpha+15)-1.$$
\end{itemize}
\end{lemma}


\medskip
\noindent \bf Proof\rm:  Let us first prove the first assertion.
Suppose that $s(t')\ne s_2(t')$ for some $t,t'$ such that
$$R(3,\alpha)\le s(t)=s_2(t)< s_2(t')\le R(3,\alpha+1),\leqno \petit{interval}$$
and $t'$ is smallest possible under (\ref{interval}). Clearly,
$t'=t+1$. Since $s(t')\ne s_2(t')$ but $s(t)=s_2(t)$, by the third
part of Proposition~\ref{sess2}, $s_2(t)+4\le s_2(t')$. This
implies that neither $s_2(t)+1$, nor $s_2(t)+2$ is a
Ramsey-number, thus $s_2(t')$ is defined from $s_2(t)$ in Case 1
(1.1, 1.2, 1.3 or 1.4) of Corollary~\ref{thm:s2}. This cannot
happen in 1.3 or in 1.4 because $s_2(t)+3<s_2(t)\le R(3,\alpha+1)$
so $s_2(t)+3$ cannot be a Ramsey number. But it cannot happen in
1.1 or in 1.2 either because there $s_2(t')\le s_2(t)+3$,
contradicting  $s_2(t)+4\le s_2(t')$ and finishing the proof.


\medskip
Now to prove the second assertion, let  $T:=\max\{t': s_2(t')<R(3,\alpha+1)\}$.
By the condition of the theorem, and the proven first part $n=s(T)=s_2(T)$. Because of Corollary~\ref{thm:s2} part  1.1,
$$s(T)\ge R(3,\alpha+1) - 2,  \,\, \hbox{and  }  T=\gap(n) =\gap(R(3,\alpha+1)).$$

Suppose for a contradiction that $s(T+i)\ne s_2(T+i)$ $(i=1,\ldots,k)$.

 By the second part of Proposition~\ref{sess2} $s(T+i)\ge s(T+i-1)+3$, so
$s(T+i)\ge s(T) + 3i\ge R(3,\alpha+1)-2 + 3i$ $(i=1,\ldots,k)$.

\medskip\noindent
{\bf Claim}: $k\le 29$

\smallskip  Indeed, otherwise $s(t+29)\ge s(t)+3\times 29\ge R(3,\alpha+1)-2 + 87= R(3,\alpha+1)+ 85.$
On the other hand, by Proposition~\ref{p:iconsecutive} (\ref{cons6})  $R(3,\alpha+1)+ 85\le R(3,\alpha+15)-1$, so
by Proposition~\ref{p:lowern}, and then applying Corollary~\ref{cor:Ramseygap}:
$$\gap_2(R(3,\alpha+1)+ 85)\ge  \left\lceil{R(3,\alpha+1)+ 85\over 2}\right\rceil - (\alpha+14)\ge
\left\lceil{R(3,\alpha+1)+1\over 2}\right\rceil +42 -
(\alpha+14)=$$
$$= \left\lceil{R(3,\alpha+1)+1\over 2}\right\rceil+42 - (\alpha+1) - 13=\gap_2(R(3,\alpha+1))+29.$$
So $s_2(t+29)\le R(3,\alpha+1)+ 85 \le s(t+29)$ and therefore there is equality throughout, proving the claim, and the theorem.
\qed

\begin{theorem}\label{thm:bounds}     For all $n,t\in\Nset$: $0\le \gap(n) - \gap_2(n) \le 2,\quad$ $0\le s_2(t) - s(t) \le 10.$
\end{theorem}

\noindent
\bf Proof\rm: Let $p<r$ two integers so that $s(p)=s_2(p)$, $s(r)=s_2(r)$, and $s(t)\ne s_2(t)$ for all $t\in\Nset$ such that $p<t<r$.
According to Lemma~\ref{thm:constant} with $\alpha:=\alpha(s(p))$,
we have $$\qquad  s(p)=R(3,\alpha+1)-1,\hbox{ or }
s(p)=R(3,\alpha+1)-2,\leqno{\petit{initial}}$$ where we can suppose
$\alpha + 1 \ge 6$ since $s(p)<R(3,5)=14$ is not possible because
then $p=3$, and the choice $t=p+1=4$ would contradict
Theorem~\ref{thm:s4}. By Lemma~\ref{thm:constant}
$$s(r)-s(p)\le 87,\qquad r-p\le 29, \qquad  \alpha(s(r))-\alpha(s(p))\le 14.$$
Moreover, for $t\in[p,r]:$  $s(t)\ge s(p)+ 3(t-p)$ that is,
$$\gap(n)\le p+ \left\lfloor{n-s(p)\over 3}\right\rfloor \leqno{\petit{ineq:gap}}$$
for any integer $n$ in the interval  $[s(p),s(r)]$, and in this same interval we are checking
$$\gap_2(n)\ge p+ \left\lfloor{n-s(p)\over 2}\right\rfloor - \beta(n-s(p)),\leqno{\petit{gap2}}$$ where $\beta: [0, 86]
\rightarrow [0,14]$ is the following function:
\begin{itemize}
\item[--] $\beta (x)=$  $0$ if $x=0$, $\beta(x)=1$  in the interval $[1,3],$
\item[--] $\beta (x)=$  $2$ in the interval $[4,7]$, $3$ in the interval $[8,11]$, $4$ in the interval $[12,17]$,
\item[--]  $5$ in $[18,22]$, $6$ in $[23,28]$, $7$ in $[29,34]$, $8$ in $[35,43]$, $9$ in $[44,48]$, $10$
in $[49,55]$,
\item[--] $11$ in $[56,62]$, $12$ in $[63,70]$, $13$ in $[71,78]$, $14$ in $[79,86].$
\end{itemize}

\bigskip We first prove  (\ref{gap2}), and then the following:

\medskip\noindent{\bf Claim } : $0\le \gap(n) - \gap_2(n) \le 2$ for all $n\in [s(p),s(r)]$.

\medskip We will be done then, since every $n \in\Nset$ belongs to such an interval by Lemma~\ref{thm:constant}.
The second assertion of the theorem also follows then: let $t\in \Nset$, and apply the first assertion to $n:=s(t)$.
Then by the proven assertion, if we are not done,  $\gap_2(n)< t=\gap(n) \le \gap_2(n)+2\le \gap_2(n+10)$, where in
the last inequality we have used the immediate consequence of Corollary~\ref{thm:s2} that $\gap_2(n+5)\ge \gap_2(n)+1$.
This means  $n < s_2(t) \le n+10=s(t)+10.$

\medskip\noindent
{\bf Proof of (\ref{gap2})}: Indeed, according to
Proposition~\ref{p:lowern}, for $n\in [s(p),s(r)]$:
$$\gap_2(n)\ge \left\lceil {n -s(p) +s(p) \over 2 }\right\rceil- \alpha(n) +\alpha -\alpha\ge
\left\lceil {s(p)\over 2 }\right\rceil -\alpha + \left\lfloor {n
-s(p)\over 2 }\right\rfloor - (\alpha(n) - \alpha) \ge$$
$$\ge p + \left\lfloor {n -s(p)\over 2 }\right\rfloor + \beta(n-s(p)),$$
where at last we applied that $\lceil {s(p)\over 2 }\rceil -\alpha=\gap_2(s_2(p))=p$ by (the first part of)
Lemma~\ref{thm:constant}; instead of the obvious estimate  $\alpha(n) - \alpha(s(p))\le \alpha (n-s(p))$
(Proposition~\ref{p:subadditivity}) we used the particular situation of the number $s(p)$ close to the
Ramsey-number $R(3,\alpha+1)$, see (\ref{initial}): the function $\beta$ provides a universal upper bound
for $\alpha(s(p)+x) - \alpha(s(p))$, {\em independently of $s(p)$}: this difference  is the number of
Ramsey-numbers in the interval $[s(p), s(p)+x]$. We have to check
$$ \alpha(s(p)+x) - \alpha(s(p)) \le \beta (x) \hbox{ for all $0\le x\le 86$.}$$

Since $\alpha\ge 4$, all the inequalities of
Proposition~\ref{p:iconsecutive} concerning $1\le k\le 14$ are
valid.  For $x=0$ the upper bound is obvious, for  $x=1$ it
follows from Proposition~\ref{p:iconsecutive} (\ref{cons1}), since
$R(3,\alpha+2)\ge R(3,\alpha+1)+3\ge s(p)+4$, for $x=2$ from
Proposition~\ref{p:iconsecutive} (\ref{cons2}), since
$R(3,\alpha+2)\ge R(3,\alpha+1)+7\ge s(p)+8,$ etc., proving
(\ref{gap2}).

\medskip\noindent
{\bf Proof of the Claim.} Of course $\gap(n)\ge \gap_2(n)$.
Combining (\ref{gap2}) and (\ref{ineq:gap}) we have
$$0\le \gap(n)-\gap_2(n) \le \left\lfloor{n-s(p)\over 3}\right\rfloor-
\left\lfloor{n-s(p)\over 2}\right\rfloor + \beta(n-s(p)),$$ which
gives our estimate by taking the maximum of the $86$ values
$x=n-s(p)$, but actually only the $14$ values
$$x=n-s(p)=1, 4, 8, 12, 18, 23, 29, 35, 44, 49, 56, 63, 71, 79,$$
matters since while  $\beta$ is constant,  the function $\gap_2(n)$
increases faster than $\gap(n)$, and the bound improves. For the
given values the differences are $$1, 1, 1, 2, 2, 2, 2, 2, 1, 2, 1,
2, 1, 1$$ in order, proving $0\le \gap(n) - \gap_2(n) \le 2$ for the
interval $[s(p),s(r)]$. \qed

  \medskip
 \noindent
{\bf Remark}:  As can be expected, the somewhat modified
computation of this proof provides the result of
Lemma~\ref{thm:constant} as well. Indeed, $\gap_2(s(p) +86)\ge p
+ 43 - 14= p +29$, that is, $s_2(p+29)\le s(p)+86$. On the other
hand, $s(p+29)\ge s(p)+ 29 \times 3=s(p)+ 87$. However,
$s_2(p+29)\ge s(p+29)$, a contradiction, proving actually $r-p\le
28$.

\bigskip
Last, we summarize the results of the two preceding theorems,
completed with the remark that the both the worst differences
between $\gap$ and $\gap_2$, $s(t)$ and $s_2(t)$ or the exception
of Theorem~\ref{thm:gap2} occur in a very small radius of
Ramsey-numbers. This can be considered as a synthesis of this
work.

\begin{theorem}\label{thm:main}  For all $n\in {\Nset}\setminus \cup _{\alpha\in \Nset}[R(3,\alpha),
R(3,\alpha)+14]:$
$\gap(n)=\gap_2(n)=\lceil n/2\rceil - \alpha(n)$, and always
  $\lceil n/2\rceil - \alpha(n)\le \gap(n) \le   \lceil n/2\rceil - \alpha(n) + 3.$

  \end{theorem}

\medskip\noindent
{\bf Proof}: The last inequality follows from the error of $2$ in
Theorem~\ref{thm:bounds} added to the additive term $1$ of
Theorem~\ref{thm:gap2}. For the first part let $\alpha\in\Nset$,
$t:=\gap_2(R(3,\alpha))$, and assume $R(3,\alpha+1) \ge
R(3,\alpha)+16$, otherwise there is nothing to prove. Then $s_2(t)<
R(3,\alpha)$ (Corollary~\ref{cor:Ramseygap}), and $s_2(t+1)\le
s_2(t)+4\le R(3,\alpha) +3$ (Corollary~\ref{thm:s2}). Set
$$I:=[s_2(t+1)+1,s_2(t+1)+12]\cap {\Nset}\subseteq [R(3,\alpha),
R(3,\alpha)+15]\subseteq [R(3,\alpha), R(3,\alpha+1) ).$$

\medskip\noindent
{\bf Claim}: If $I$  does not contain any Ramsey-number, then
there exists $t'\in \Nset$:
$$s(t')=s_2(t')\in I.$$

\medskip Indeed,  by the condition $\alpha$ is constant on $I$, so by Theorem~\ref{thm:gap2},
$s_2(t+7)= s_2(t+1)+12\le R(3,\alpha)+15$.  On the other hand, by
Proposition~\ref{sess2} we have $s(t+7)\ge s(t+3)+12$. If the Claim
is not true, the equality does not hold here, whence
$s_2(t+1)>s(t+3)$. This means that defining  $n=s(t+3)$, we have
$\gap_2(n)\le t$ and $\gap(n)=t+3$, contradicting
Theorem~\ref{thm:bounds} and proving the claim.

\medskip Now by Lemma~\ref{thm:constant}, for the $t'$ provided by the Claim and
for any $n \in [s_2(t'), R(3,\alpha+1)]$ we have
$\gap_2(n)=\gap(n)$. According to the Claim, $s(t')\le
R(3,\alpha)+15$ finishing the proof.
\qed

\section{Graphs with small gap}\label{sub:start}

In this section we explore the smallest gap-extremal graphs and for
small orders we show the graphs of maximum gap.  Graphs on at most $4$
vertices are perfect, so $s(1)\ge 5$, and the only $1$-extremal graph
is $C_5$.

We will need the following lemma of merely technical use.  A graph $G$
is \emph{clique-Helly} if its inclusion-wise maximal cliques (viewed as set of vertices) have
the Helly property: if a collection of maximal cliques of $G$ pairwise
intersect, then they have a common vertex.  A triangular claw is a
graph $T_6$ on $6$ vertices, and $9$ edges consisting of a triangle
$\Delta\subseteq V(T_6)$ and a $3$-stable set $S\subseteq V(T_6)$,
$V(T_6)=\Delta\cup S$ so that every vertex of $S$ is joined to a
different pair of vertices of $\Delta$.  This graph is not
clique-Helly, and as shown below, it is in a sense the basic example of
a non-clique-Helly graph.  We omit the simple proof of the following Lemma:

\begin{lemma}[See \cite{LinS:07}]\label{l:triangularclaw}
  If a graph $G$ does not contain a triangular claw as an induced subgraph then
  it is clique-Helly.
\end{lemma}

\begin{theorem}\label{thm:s210}
  The graph $2C_5$ is gap-extremal, in particular, $s(2)=s_2(2)=10$
  and the only $2$-extremal graph is $2C_5$. Therefore the graphs
  consisting of a $C_5$ and an arbitrary graph on $\{1\}$, $\{1,2\}$,
  $\{1,2,3\}$, $\{1,2,3,4\}$ have maximum gap for $n=6,7,8,9$
  respectively.  In addition
\begin{itemize}
\item[(1)] for $n=6$ this is the unique graph of maximum gap, and
it is stable-gap-optimal. \item[(2)] for $n=7$ the gap of $C_7$
and $\bar C_7$ is maximum, as well as that of $R-v$ where $R$ is a
$(3,4)$-Ramsey graph and $v\in V(R)$. The latter graphs are
stable-gap-optimal.
\item[(3)] for $n=8$  the only
stable gap-optimal graphs are the $(3,4)$-Ramsey graphs.
\item[(4)] for $n= 9$  a graph $G$ on $n$ vertices is
stable gap-optimal if and only if it is triangle-free and
$\alpha(G)=4$.
\end{itemize}
\end{theorem}

\medskip\noindent{\bf Proof}: We first prove (1) and (2). By Proposition~\ref{l:jumpk} $s(2)\ge s(1)+2=7$,
so $\gap(6)=\gap(7)=1$, and (2) immediately follows.  A graph $G$ of maximum gap on $6$ vertices is imperfect,
so it contains $C_5$ as induced subgraph. The vertex $v$ not contained in this $C_5$ is an isolated vertex,
since otherwise the edge $vu$ and the matching of $C_5-u$ is a clique cover with $3$ edges, whence $\gap(G)=0$,
a contradiction which proves (1).

Suppose that $G$ is a $2$-extremal  graph. Since $\gap (2C_5)=2,$
we have $n:=|V(G)|\le 10$. The only thing we have to prove now is
$G=2C_5$, since then $\gap(8)=\gap(9)=1$ follow and (3) and (4)
can be readily checked: by Proposition~\ref{p:firstramsey}
$\alpha(8)=3$, $\alpha(9)=4$, so for  any triangle-free graph $G$
on $8$ vertices with $\alpha(G)=3$ we have    $\gap(G)\ge \lceil
n/2 \rceil  - 3=1$, and for  any triangle-free graph $G$ on $9$
vertices with $\alpha(8)=4$ we have    $\gap(G)\ge \lceil n/2
\rceil  - 3=1$. It follows that their gap is maximum, and on $8$
vertices these are exactly the $(3,4)$-Ramsey graphs. Conversely,
stable-gap-optimal graphs are triangle-free and their stability
number is as claimed by definition, so the assertion follows from
the proven part.

Suppose now for a contradiction that $G\ne 2C_5$. Let
$\alpha:=\alpha(G)$, $\omega:=\omega(G)$, $\theta:=\theta(G)$.

\medskip
\noindent {\bf Claim 1}: If $K$ is a clique of $G$, then $G-K$ has
at least $7$ vertices.

By Proposition~\ref{l:removeK} $\gap(G-K)=1$, so it has at least $5$
vertices.  If it has exactly $5$ vertices, then it is a $C_5$. Then
$\theta(G)\le 4$, so $\alpha(G)\le \theta(G)-2\le 2$, and the equality
holds everywhere. Pick a vertex $v$ of this $C_5$. Then $N(v)$ is the
union of a stable set and a clique, so it does not contain a $C_5$,
$C_7$ or $\bar C_7$ (it is \emph{split graph}), so $N(v)$ induces a
perfect graph, and we conclude $\gap(G)\le 1$ by
Proposition~\ref{l:alpha2}.  If $G-K$ has $6$ vertices, then by (1)
$G-K$ has an isolated vertex $v$, whence $N(v)$ is simplicial in $G$,
contradicting Proposition~\ref{simplicial}.

 \medskip\noindent {\bf Claim 2}: $\alpha=\omega=3$, $\theta=5$, $n=10$.

Apply Claim~1 to an arbitrary clique $K$. Since $n\le 10$, we get
$|K|\le 3$. If there exists a clique $K$ for which equality holds,
we have $n=10$, $\omega=3$.

If $\omega\le 2$, then by Proposition~\ref{l:fcrit} every component
of $G$ is factor-critical, that is odd, and at least two of them are
imperfect: $G=2C_5$. So $\omega=3$ and $n=10$.

Now by Proposition~\ref{p:firstramsey}, $R(4,3)=R(3,4)=9$, so
since $\omega=3$,  $\alpha\ge 3$. But $\alpha\ge 4$ is not
possible, because then by  Proposition~\ref{l:removeK}
$\alpha(G-K)\ge 4,\, \gap(G-K)=1,\, \theta(G-K)\ge 5.$ Since $G-K$
has $7$ vertices but is neither $C_7$ nor $\bar C_7$, it contains
a $C_5$, and the two  vertices that are not in this $C_5$ are
isolated ones because of $\theta(G-K)\ge 5$. If $v$ is one of
them, then again, it is a simplicial vertex in $G$, contradicting
Proposition~\ref{simplicial}, and finishing the proof of the
claim.

\medskip\noindent {\bf Claim 3}: $G$ contains two disjoint triangles.

\smallskip Because of $\theta(G-v)=4$, we have $\omega(G-v) = 3$ for
all $v\in V$.  If $G$ does not contain 2 disjoint triangles, then the
triangles of $G$ pairwise intersect, so either $G$ is clique-Helly and
they all intersect, a contradiction to $\omega(G-v) = 3$, or by
Lemma~\ref{l:triangularclaw}, $G$ contains a triangular claw $\Delta
\cup S$ where $S=\{s_1,s_2,s_3\}\subseteq V(G)$ is a stable set, and
$\Delta=\{t_1,t_2,t_3\}\subseteq V(G)$ is a triangle, and $s_i$ is
adjacent to $T\setminus \{t_i\}$ $(i=1,2,3)$.  Note that $\Delta \cup
S$ may be assumed to be induced because adding an edge to it yields
either a $K_4$ or two disjoint triangles.

We may assume that $G - \{t_1, t_2\}$ is triangle-free because else,
there are two disjoint triangles.  Since $\alpha(G - \{t_1, t_2\}) =
3$, $G - \{t_1, t_2\}$ must be one of $W_8$, $W_{81}$, $W_{82}$
(Proposition~\ref{p:firstramsey}).  So, $G - \{t_1, t_2\}$ has
a cycle  $w_1 \dots w_8w_1$, and the only other edges are among
$w_iw_{i+5}$, $i=1, \dots, 4$.  We suppose up to symmetry  $t_3=w_1$.
We consider now two cases.

{\bf Case 1}, $t_1$ is not adjacent to $w_2$ and $w_8$.  Because of
the triangular claw, $w_1$ and $t_1$ have a common neighbor that must
be $w_5$.  Also $t_2$ and $w_1$ must have a common neighbor, that
cannot be $w_5$ because $\omega=3$, so it is $w_2$ or $w_8$, say $w_2$
up to symmetry.  Now, we may assume $t_2w_3, t_1w_4, t_1w_6 \notin
E(G)$ because otherwise there are two disjoint triangles.  So, the
common neighbor $s_3$ of $t_1t_2$ must be $w_7$ and we may assume
$t_2w_6, t_2w_8 \notin E(G)$ because otherwise there are two disjoint
triangles.  Hence, $\{t_2, w_3, w_6, w_8\}$ is a stable set, a
contradiction.

{\bf Case 2}, $t_1$ has at least one neighbor among $w_2$ and $w_8$.
Symmetrically, we may assume that $t_2$ also has at least one neighbor
among $w_2$ and $w_8$.  Since $\omega = 3$, we may assume $t_1w_8,
t_2w_2 \in E(G)$ and $t_1w_2, t_2w_8 \notin E(G)$.  Now, we may assume
$t_1w_7, t_2w_3 \notin E(G)$ because otherwise there are two disjoint
triangles.  Hence, $\{t_1, w_2, w_7, w_4\}$ is a stable set unless
$t_1w_4 \in E(G)$, so $t_1w_4 \in E(G)$ and symmetrically, $t_2w_6 \in
E(G)$.  Now, $t_2w_5\notin E(G)$ because else there are two disjoint
triangles.  Hence, $\{t_2, w_3, w_5, w_8\}$ is a stable set, a
contradiction.  This proves the claim.

\medskip So, $G$ contains two vertex-disjoint triangles,
$T_1=\{a_1,a_2,a_3\}, T_2=\{b_1,b_2,b_3\}$.  If the remaining four
vertices contain a triangle or two independent edges, we have
$\theta(G)\le 4$, a contradiction.  Therefore three of these vertices
form an independent set $C=\{c_1,c_2,c_3\}$ and we have the following
cases according to the adjacencies of the last vertex~$d$ (which has a
neighbor among $c_1, c_2, c_3$ because $\alpha(G) = 3$).

\bf Case 1, \rm $dc_i\in E(G)$ for $i=1,2,3$.  Each vertex of $T_1$
must have a neighbor in $C$ because $\alpha(G) = 3$.  If $a_1c_1,
a_2c_1 \in E(G)$ then we must have $a_3c_2 \in E(G)$ or $a_3c_3\in
E(G)$ because there is no $K_4$. But then, we can cover $G$ with two
triangles and two edges.  So we proved that no two vertices in $T_1$
can have a common neighbor in $C$.  Hence, we may assume that the only
edges between $T_1$ (and similarly $T_2$) and $C$ are $c_ia_i$ (and
similarly $c_ib_i$), $i=1,2,3$.  Using that $\alpha(G)=3$, it follows
that $a_ib_i\in E(G)$ and now $a_i,b_i,c_i$ for $i=1,2,3$ give three
disjoint triangles showing that $\theta(G)\le 4$, a contradiction.

\bf Case 2, \rm $dc_3\in E(G),dc_1,dc_2\notin E(G)$.  Suppose first
that every vertex of $T_1$ has a neighbor in $\{c_1, c_2\}$.  Since
there is no $K_4$ we may assume $a_1c_1, a_2c_2, a_3c_2 \in E(G)$, so
we can cover $G$ with two triangles and two edges, a contradiction.
So there must be a vertex in $T_1$ with no neighbor in $\{c_2, c_1\}$,
say $a_1$, and by the same argument a similar vertex in $T_2$, say
$b_1$.  Using five times that $\alpha(G)=3$, we get that
$a_1c_3,b_1c_3,a_1b_1,da_1,db_1\in E(G)$, a contradiction because
$\{a_1,b_1,c_3,d\}$ is a clique.

\bf Case 3, \rm $dc_2,dc_3\in E(G),dc_1\notin E(G)$.  We claim that
$c_1$ is nonadjacent to at least two vertices of both $T_1,T_2$.  If
not, say $c_1$ is adjacent to $a_2, a_3$, then $c_2a_1,c_3a_1 \notin
E(G)$ otherwise we have a cover with two triangles and two
edges. Depending on $c_1a_1\in E(G)$ or not, we have either a clique
or an independent set of size four, a contradiction that proves the
claim.  Therefore, w.l.o.g.\ $c_1$ is non-adjacent to $a_2, a_3, b_2,
b_3$.  If $c_1a_1\notin E(G)$ or $c_1b_1\notin E(G)$ or $a_1b_1 \in
E(G)$ then $c_1$ is a simplicial vertex, a contradiction.  Thus
$c_1a_1, c_1b_1\in E(G), a_1b_1\notin E(G)$.

  Next we note that each of $a_2, a_3$ must have a neighbor in $\{c_2,
  c_3\}$, else there is an $S_4$.  But $a_2, a_3$ may not have a
  common neighbor in $\{c_2, c_3\}$ because then there is a cover with
  two triangles and two edges.  Hence w.l.o.g.\ the only edges between
  $T_1$ and $C$ are $c_1a_1, c_2a_2, c_3a_3$.  Similarly, the only
  edges between $T_2$ and $C$ are $c_1b_1, c_2b_2, c_3b_3$.

  Now $\alpha(G)=3$ implies $a_2b_2,a_3b_3\in E(G)$.  Moreover $da_2,
  da_3, db_2, db_3 \notin E(G)$ otherwise there is a clique cover with
  two triangles and two edges. Then $a_2b_3, a_3b_2 \in E(G)$ for
  otherwise $a_2, b_3, c_1, d$ or $a_3, b_2, c_1, d$ would form an
  independent set.  But now have the final contradiction since $a_2,
  a_3, b_2, b_3$ span a clique.
\qed

 \bigskip To slightly shorten the proof, one could use Chv\'atal's
 \cite{chvatal:74} theorem stating that the \emph{Gr\"otzsch graph}
 (the fourth in Mycielski's well-known construction \cite{M}, being
 the ``Mycielskian" of $C_5$ which is the third) is {\em the only
   triangle-free graph on at most $11$ vertices with chromatic number
   at least $4$.}  The complement of the Gr\"otzsch graph is therefore
 the only graph on at most $11$ vertices with $\alpha\le 2$ and
 $\theta\ge 4$.  Also the following lemma could be used.  For a proof,
 see Lemma 1.16 in \cite{nicolas:hdr}.

\begin{lemma}
  \label{threecases}
  If $G$ is a graph on at least 10 vertices then either $G$ contains
  a clique or a stable set on four vertices, or $G$ contains two
  disjoint triangles.
\end{lemma}

\begin{theorem}\label{thm:s3}
The graph $R_{13}$ is gap-extremal, in particular,
$s(3)=s_2(3)=13$, and   the only $3$-extremal graph is $R_{13}$.
Any triangle-free graph $G$ on $11$ or $12$ vertices and
$\alpha(G)\le 4$ is stable-gap-optimal and connected.
\end{theorem}

\medskip\noindent
\bf Proof\rm :  Suppose that $G$ is a $3$-extremal  graph,
$\alpha:=\alpha(G)$, $\omega:=\omega(G)$, $\theta:=\theta(G)$. Since
$R_{13}$  is triangle-free, $\theta(R_{13})=\zeta(R_{13})=7$, and
$\alpha(R_{13})=4$ (it is a $(3,5)$-Ramsey graph).  So
$\gap(R_{13})=3$, and therefore $n:=|V(G)|\le 13$. We have to prove
$G=R_{13}$. If $\omega= 2$ this is true since then by
Proposition~\ref{l:fcrit} $G$ is factor-critical, $\theta(G)=7$, so
$\alpha(G)=4$. Therefore $G$ is a $(3,5)$-Ramsey graph, and by
Proposition~\ref{p:firstramsey} $G=R_{13}$.
 So suppose $\omega\ge 3$.

\medskip\noindent{\bf Claim 1}: $n=13,\,\omega=3,\,\alpha=4,\, \theta=7,$ and for every triangle
$T$, $G-T$ is a $2C_5$.

\smallskip
If $K$ is an arbitrary clique, $\gap (G-K) = 2$, so by
Theorem~\ref{thm:s210}, $G-K$ is of order at least $10$, whence
$|K|\le 3$, and therefore $\omega=3$. If $T$ is a triangle, $n\le
13$ implies that $G-T$ is of order at most $10$. So $G-T$ is of
order $10$ and gap $2$, and $n=13$. By
Proposition~\ref{l:removeK}, $\gap(G-T)=2$, and since $G-T$ has
$10$ vertices, the unicity in Theorem~\ref{thm:s210} states that
it is $2C_5$.

Now by the equalities of Proposition~\ref{l:removeK} concerning
gap-critical graphs, $\alpha(G) =\alpha(G -
Q)=\alpha(2C_5)=4=\theta-3$, finishing the proof of the claim.

\medskip\noindent{\bf Claim 2}: Let $T$ be a triangle, and
$C,D\subseteq V(G)$ be the two $C_5$ components of $G-T$. Then for
every $t\in T$ either $\alpha(\{t\}\cup C)=2$ or $\alpha(\{t\}\cup
D)=2$.

\smallskip\noindent Indeed, if there exists $t\in T$ so that both are
$3$, then there exists $c_1, c_2\in C$, and $d_1, d_2\in D$ so that
$t, c_1, c_2, d_1, d_2$ form a stable set in $G$, contradicting Claim
2.

\medskip\noindent So suppose $t\in T$, $\alpha(\{t\}\cup C)=2$. Then
$C\sm N(t)$ is the subset of an edge of $C$, and therefore $t$ forms a
triangle $T_1$ and $T_2$ with two different edges of $C$. But this is
impossible, because by Claim~1 both $G-T_1$ and $G-T_2$ are $2C_5$
graphs, however, $(C-T_1)\cup T\setminus \{t\}\ne (C-T_2)\cup
T\setminus \{t\}$, because $T_1\ne T_2$.

The remaining additional claim follows now from
Proposition~\ref{p:lower}: if $G$ is a triangle-free graph on $11$ or
$12$ vertices and $\alpha(G)=\alpha(n)=4$, then $\gap(G)\ge 6 - 4=2$,
so the equality holds and $\gap (G)$ is maximum. Moreover $G$ is
connected since $R(3,2)=3$, $R(3,3)=6$, $R(3,4)=9$ imply that two
vertex-disjoint graphs with stability numbers $2$ and $2$ or $1$ and
$3$ have at most 12 vertices. \qed

\begin{theorem}\label{thm:s4}
The $(3,6)$-Ramsey graphs are $4$-extremal, in particular
$s(4)=s_2(4)=17$. A graph is $4$-extremal and triangle-free if and
only if it is a $(3,6)$-Ramsey graph; for all other (possibly
non-existing) $4$-extremal graphs $G$, $\alpha(G)=4$, and
$\theta(G)=8$.

\end{theorem}

\medskip\noindent
\bf Proof. \rm  Let $G$ be $4$-extremal. According to
Proposition~\ref{p:lower} the gap of  $(3,6)$-Ramsey graphs on
$17$ vertices is at least $9-5=4$. So $n:=|V(G)|\le 17$.


 Since $s_2(4)=17$ from Corollary
\ref{values}, we may assume $\omega:=\omega(G)\ge 3$. Then by
Proposition~\ref{l:jumpk} $s(4)\ge s(3)+ 3=16$ (see
Theorem~\ref{thm:s3}), and $s(4)\le s_2(4)=17$. The statement
$s(4)=17$ follows now from the next claim.

\medskip\noindent{\bf Claim}: For any clique $K$, $G-K$ is of order at least $13$, $\omega\le 4$, and $|V(G)|=17$.

\smallskip
Indeed, by Proposition~\ref{l:removeK}, $\gap(G-K)=3$. So $G-K$ is
of order at least $13$, so $|K|\le 4$.  Suppose $n=16$. Then
$\omega=3$ and for any triangle $K$, $G-K$ is of order exactly
$13$ of gap $3$, so it is a $(3,5)$-Ramsey graph, in particular it
is triangle free. Consequently there are no two disjoint triangles
in $G$, and $\alpha(G)=\alpha(G-K)=\alpha(R_{13})=4$.

On the other hand $n-2= 14=R(3,5)$, so for all $u,v\in V(G)$,
$\omega(G-\{u,v\})=3$.

So, by Lemma~\ref{l:triangularclaw}, $G$ is clique-Helly or has a
triangular claw.  In the first case, since there are no two disjoint
triangles, the triangles pairwise intersect, so they intersect, a
contradiction to $\omega(G-\{u,v\})=3$.  Hence, there is a triangular
claw $\{t_1, t_2, t_3, s_1, s_2, s_3\}$ (our usual notation).  Since
$\omega(G-\{t_1, t_2\})=3$, $G-\{t_1, t_2\}$ contains a triangle,
hence $G$ contains two disjoint triangles.  This contradiction
finishes the proof of the Claim.

\smallskip

Let $K$ be an $\omega$-clique of $G$. By
Proposition~\ref{p:gallaicol}, $\theta(G-K)\le 7$, so by
Proposition~\ref{l:removeK} $\theta(G)\le 8$, and since
$\gap_2(G)=4$: $\alpha(G)\le 4$. The strict inequality here, that
is, $\alpha \le 3$ would imply either $\omega\le 3$ and then
applying $R(4,4)=18$ (Proposition~\ref{p:firstramsey}) we get that
$G$ is a $(4,4)$-Ramsey-graph; or  by Claim, $\omega= 4$, and $G-K$
is of gap $3$ and order $13$, so isomorphic to $R_{13}$. In the
former case we see that $\theta=6$, implying $\alpha=2$, but
$(4,4)$-Ramsey graphs have $\alpha=3$, a contradiction; in the
latter case $\alpha(G)=\alpha(G-K)=4$ is proved, finishing the proof
of the theorem. \qed

\bigskip
Surprisingly, the next case we can treat is $s(10)$:

\begin{lemma}\label{lem:s10}
$s(t)\ge s_2(4) + 3(t-4)$ for $t=5,\ldots, 10.$
\end{lemma}

\noindent \bf Proof\rm: Note that $s_2(i)-s_2(i-1)$ for the six
values $i=5,\ldots, 10$ is equal to $4, 4, 4, 2, 2, 2$, that is,
$3$ in average.

If the statement does not hold let $t_0$ be the smallest value for
which this inequality is violated. Then $$s(t_0)< s_2(4) +
3(t_0-4) \le s_2(t_0).$$ Clearly, $s(t_0)-s(t_0-1)=2$ since if
not, according to Proposition~\ref{sess2} $s(t_0)-s(t_0-1)\ge 3$
so we could have chosen $t_0-1$ or a smaller value instead of
$t_0$. Therefore any $t_0$-extremal graph is triangle free, in
contradiction with $s(t_0)< s_2(t_0)$.\qed

\medskip
Using Lemma \ref{lem:s10} for a lower bound and Corollary
\ref{values} as upper bound, $s(5)\in \{20,21\}$, $s(6)\in
\{23,24,25\}$, $s(7)\in \{26,27,28\}$, $s(8)\in \{29,30,31\}$,
$s(9)\in \{32,33\}$, and $s(10)=35$.

\begin{corollary}\label{thm:s10}
We have  $s(10)=35$, the $(3,9)$-Ramsey graphs are all
$10$-extremal, and all other $10$-extremal graphs contain a
triangle.
\end{corollary}

\noindent \bf Proof\rm:  $s(10)\le 35,$ since by
Proposition~\ref{p:lowern} $\gap_2(35)\ge \lceil 35/2 \rceil -
\alpha (35)=18 - 8=10$, so $s_2(10)\le 35)$. Substituting
$s(4)=17$ (Theorem~\ref{thm:s4}) and $t=10$  into
Lemma~\ref{lem:s10} we get $s(10)\ge 35$. \qed

\medskip
\noindent {\bf Acknowledgment}: Thanks to Zoli F\"uredi for calling
our attention to \cite{BFJ} (and its references).

\end{document}